\documentclass[12pt]{article}
\usepackage{amssymb, amsmath, amsthm, amscd}
\usepackage[pdftex]{graphics}
\usepackage[utf8]{inputenc}
\usepackage[all,cmtip]{xy}
\usepackage{enumitem}
\usepackage{setspace}
\usepackage{mathdots}

\usepackage[mathscr]{euscript}
\usepackage[colorlinks=true]{hyperref}
\hypersetup{colorlinks   = true}
\hypersetup{linkcolor=blue}
\usepackage{tikz}
\begin{document}

\mathchardef\mhyphen="2D
\newtheorem{The}{Theorem}[section]
\newtheorem{Lem}[The]{Lemma}
\newtheorem{Prop}[The]{Proposition}
\newtheorem{Cor}[The]{Corollary}
\newtheorem{Rem}[The]{Remark}
\newtheorem{Obs}[The]{Observation}
\newtheorem{SConj}[The]{Standard Conjecture}
\newtheorem{Titre}[The]{\!\!\!\! }
\newtheorem{Conj}[The]{Conjecture}
\newtheorem{Question}[The]{Question}
\newtheorem{Prob}[The]{Problem}
\newtheorem{Def}[The]{Definition}
\newtheorem{Not}[The]{Notation}
\newtheorem{Claim}[The]{Claim}
\newtheorem{Conc}[The]{Conclusion}
\newtheorem{Ex}[The]{Example}
\newtheorem{Fact}[The]{Fact}
\newtheorem{Formula}[The]{Formula}
\newtheorem{Formulae}[The]{Formulae}
\newtheorem{The-Def}[The]{Theorem and Definition}
\newtheorem{Prop-Def}[The]{Proposition and Definition}
\newtheorem{Lem-Def}[The]{Lemma and Definition}
\newtheorem{Cor-Def}[The]{Corollary and Definition}
\newtheorem{Conc-Def}[The]{Conclusion and Definition}
\newtheorem{Terminology}[The]{Note on terminology}
\newcommand{\C}{\mathbb{C}}
\newcommand{\R}{\mathbb{R}}
\newcommand{\N}{\mathbb{N}}
\newcommand{\Z}{\mathbb{Z}}
\newcommand{\Q}{\mathbb{Q}}
\newcommand{\Proj}{\mathbb{P}}
\newcommand{\Rc}{\mathcal{R}}
\newcommand{\Oc}{\mathcal{O}}
\newcommand{\Vc}{\mathcal{V}}
\newcommand{\Id}{\operatorname{Id}}
\newcommand{\pr}{\operatorname{pr}}
\newcommand{\rk}{\operatorname{rk}}
\newcommand{\del}{\partial}
\newcommand{\delbar}{\bar{\partial}}
\newcommand{\Cdot}{{\raisebox{-0.7ex}[0pt][0pt]{\scalebox{2.0}{$\cdot$}}}}
\newcommand\nilm{\Gamma\backslash G}
\newcommand\frg{{\mathfrak g}}
\newcommand{\fg}{\mathfrak g}
\newcommand{\Oh}{\mathcal{O}}
\newcommand{\Kur}{\operatorname{Kur}}
\newcommand\gc{\frg_\mathbb{C}}
\newcommand\jonas[1]{{\textcolor{green}{#1}}}
\newcommand\luis[1]{{\textcolor{red}{#1}}}
\newcommand\dan[1]{{\textcolor{blue}{#1}}}

\begin{center}

{\Large\bf Some Properties of Balanced Hyperbolic Compact Complex Manifolds}

\end{center}

\begin{center}

{\large Samir Marouani and Dan Popovici}

\end{center}

\vspace{1ex}

\noindent{\small{\bf Abstract.} We prove several vanishing theorems for the cohomology of balanced hyperbolic manifolds that we introduced in our previous work and for the $L^2$ harmonic spaces on the universal cover of these manifolds. Other results include a Hard Lefschetz-type theorem for certain compact complex balanced manifolds and the non-existence of certain $L^1$ currents on the universal covering space of a balanced hyperbolic manifold.}

\vspace{2ex}

\section{Introduction}\label{section:Introduction} In this paper, we continue the study of compact complex {\it balanced hyperbolic} manifolds that we introduced very recently in [MP21] as generalisations in the possibly non-projective and even non-K\"ahler context of the classical notion of {\it K\"ahler hyperbolic} (in the sense of Gromov) manifolds. Recall that every K\"ahler hyperbolic manifold is also {\it Kobayashi/Brody hyperbolic}.

%\vspace{2ex}

%Let $X$ be a {\it compact} complex manifold. Recall that $X$ is said to be {\it Brody hyperbolic} ([Bro78, Theorem 4.1]) if there are no entire curves in $X$, namely there exists no non-constant holomorphic map $f:\C\longrightarrow X$. Since $X$ is compact, this is known to be equivalent to $X$ being {\it Kobayashi hyperbolic} in the sense that the Kobayashi pseudo-distance on $X$ is actually a distance. (See e.g. [Kob70].) In [MP21, Definition 2.5], we introduced the following $1$-codimensional analogue of this: \\

%\noindent {\it Let $n\geq 2$ be an integer. An $n$-dimensional compact complex manifold $X$ is said to be {\bf divisorially hyperbolic} if there exists no holomorphic map $f:\C^{n-1}\longrightarrow X$ such that $f$ is non-degenerate at some point $x\in\C^{n-1}$ and has subexponential growth.} \\

%What we mean by $f$ having a {\it subexponential growth} is spelt out in Definition 2.3 of [MP21]. It refers to the growth of the volume of the ball $B_r\subset\C^{n-1}$ of radius $r>0$ centred at $0\in\C^{n-1}$, with respect to the degenerate metric $f^\star\omega$ on $\C^{n-1}$ that is the pullback under $f$ of an arbitrary Hermitian metric $\omega$ on $X$, as $r$ tends to $+\infty$. Since we do not deal with {\it divisorially hyperbolic} manifolds in this paper, we do not repeat here the details of Definition 2.3 of [MP21].

\vspace{2ex}

On the other hand, recall that a Hermitian metric on a complex manifold $X$ identifies with a $C^\infty$ positive definite $(1,\,1)$-form $\omega$ on $X$. If we put $\mbox{dim}_\C X=n\geq 2$, a Hermitian metric $\omega$ is said to be {\bf balanced} (see [Gau77] where these metrics were introduced under the name of semi-K\"ahler and [Mic83] where they were given this name) if $d\omega^{n-1}=0$. Moreover, $\omega$ is said to be {\bf degenerate balanced} (see [Pop15] for the name) if $\omega^{n-1}$ is $d$-exact. Unlike in the K\"ahler setting, where no $d$-exact Hermitian metric $\omega$ can exist on a {\it compact} complex manifold, compact complex manifolds carrying degenerate balanced metrics do exist. These manifolds include:

\vspace{1ex}

(i)\, the connected sums $X=\sharp_k(S^3\times S^3)$, with $k\geq 2$, endowed with the Friedman-Lu-Tian complex structure constructed via conifold transitions ([Fri89], [LT93], [FLY12]);

\vspace{1ex}

(ii)\, the quotients $X=G/\Gamma$ of any {\bf semi-simple} complex Lie group G by a lattice $\Gamma\subset G$ ([Yac98]).

\vspace{1ex}

In [MP21], we generalised the notion of degenerate balanced compact complex manifolds starting from the observation that this is a kind of hyperbolicity property.

Throughout the text, $\pi_X:\widetilde{X}\longrightarrow X$ will stand for the universal cover of $X$. If $\omega$ is a Hermitian metric on $X$, we let $\widetilde\omega=\pi_X^\star\omega$ be the Hermitian metric on $\widetilde{X}$ that is the lift of $\omega$. According to [Gro91], a $C^\infty$ $k$-form $\alpha$ on $X$ is said to be $\widetilde{d}(\mbox{bounded})$ with respect to $\omega$ if $\pi_X^\star\alpha = d\beta$ on $\widetilde{X}$ for some $C^\infty$ $(k-1)$-form $\beta$ on $\widetilde{X}$ that is bounded w.r.t. $\widetilde\omega$.

Now, recall that a compact complex manifold $X$ is said to be {\it K\"ahler hyperbolic} in the sense of Gromov (see [Gro91]) if there exists a K\"ahler metric $\omega$ on $X$ (i.e. a Hermitian metric $\omega$ with $d\omega=0$) such that $\omega$ is $\widetilde{d}(\mbox{bounded})$ with respect to itself. In [MP21, Definition 2.1], we introduced the following 1-codimensional analogue of this:\\

\noindent {\it An $n$-dimensional compact complex manifold $X$ is said to be {\bf balanced hyperbolic} if there exists a balanced metric $\omega$ on $X$ such that  $\omega^{n-1}$ is $\widetilde{d}(\mbox{bounded})$ with respect to $\omega$.} \\

\noindent Any such metric $\omega$ is called a {\it balanced hyperbolic} metric. 
\vspace{2ex}

The implications among these notions are summed up in the following diagram. (See [MP21] for relations with other hyperbolicity notions.)

\vspace{3ex}

\hspace{12ex} $\begin{array}{lll} X \hspace{1ex} \mbox{is {\bf K\"ahler hyperbolic}} & \implies & X \hspace{1ex} \mbox{is {\bf balanced hyperbolic}} \\
  &  & \rotatebox{90}{$\implies$} \\
     &  & X \hspace{1ex} \mbox{is {\bf degenerate balanced}} \end{array}$

\vspace{3ex}

We now outline the specific properties of these classes of manifolds that we prove in this paper.

\vspace{2ex}

{\bf (I) Case of balanced and degenerate balanced manifolds}

\vspace{1ex}

In the first part of the paper, we obtain some general results on compact complex manifolds carrying {\it balanced} metrics (and, in some cases, results on {\it Gauduchon} metrics) and then we use them to infer vanishing results for {\it degenerate balanced} manifolds. See $\S.$\ref{subsection:deg-bal_background} for a reminder of the terminology used in what follows.

\vspace{2ex}

{\bf (a)}\, Our first main result, obtained as a consequence of the computation in Lemma \ref{Lem:1-forms_Delta-harm}, is a {\bf Hard Lefschetz Isomorphism} between the De Rham cohomologies of {\bf degrees 1} and {\bf 2n-1} that holds on any compact complex {\it balanced} manifold satisfying a mild $\partial\bar\partial$-type condition.

For any Hermitian metric $\omega$ on an $n$-dimensional complex manifold, we will use throughout the paper the notation $\omega_p:=\omega^p/p!$ for any integer $p$ between $2$ and $n$.

\begin{The}\label{The:Hard-Lefschetz_bal-ddbar} Let $X$ be a compact complex manifold with $\mbox{dim}_\C X=n$.

  \vspace{1ex}

  (i)\, If $\omega$ is a {\bf balanced} metric on $X$, the linear map: \begin{equation}\label{eqn:Hard-Lefschetz_map}\{\omega_{n-1}\}_{DR}\wedge\cdot : H^1_{DR}(X,\,\C)\longrightarrow H^{2n-1}_{DR}(X,\,\C), \hspace{3ex} \{u\}_{DR}\longmapsto\{\omega_{n-1}\wedge u\}_{DR},\end{equation} is well defined and depends only on the cohomology class $\{\omega_{n-1}\}_{DR}\in H^{2n-2}_{DR}(X,\,\C)$.

  \vspace{1ex}

(ii)\, If, moreover, $X$ has the following additional property: for every form $v\in C^\infty_{1,\,1}(X,\,\C)$ such that $dv=0$, the following implication holds: \begin{equation}\label{eqn:Hard-Lefschetz_bal-ddbar}v\in\mbox{Im}\,\partial\implies v\in\mbox{Im}\,(\partial\bar\partial),\end{equation} the map (\ref{eqn:Hard-Lefschetz_map}) is an {\bf isomorphism}.

\end{The}

As a consequence of this discussion, we obtain the following vanishing properties for the cohomology of degenerate balanced manifolds.

\begin{Prop}\label{Prop:BC_10-01_b1_deg-bal} Let $X$ be a compact {\bf degenerate balanced} manifold.

 \vspace{1ex} 

 (i)\, The Bott-Chern cohomology groups of types $(1,\,0)$ and $(0,\,1)$ of $X$ vanish: $H_{BC}^{1,\,0}(X,\,\C)=0$ and $H_{BC}^{0,\,1}(X,\,\C)=0$.

 \vspace{1ex} 

(ii)\, If, moreover, $X$ satisfies hypothesis (\ref{eqn:Hard-Lefschetz_bal-ddbar}), its De Rham cohomology group of degree $1$ vanishes: $H^1_{DR}(X,\,\C)=0$. 

\end{Prop}

Note that degenerate balanced manifolds that satisfy hypothesis (\ref{eqn:Hard-Lefschetz_bal-ddbar}) do exist. Indeed, Friedman showed in [Fri19] that the manifolds $X=\sharp_k(S^3\times S^3)$, with $k\geq 2$, endowed with the Friedman-Lu-Tian complex structure constructed via conifold transitions ([Fri89], [LT93], [FLY12]) are even $\partial\bar\partial$-manifolds.

\vspace{1ex}

{\bf (b)}\, Our study of the cohomology of {\bf degree 2} in this setting centres on seeking out possible positivity properties of balanced hyperbolic manifolds. As an alternative to Question 1.4. in [MP21], wondering about possible positivity properties, in the senses defined therein, of the canonical bundle $K_X$ of any balanced hyperbolic manifold $X$, we concentrate this time on whether there are ``many'' (in a sense to be determined) closed positive currents $T$ of bidegree $(1,\,1)$ on such a manifold.

The starting point of this investigation is Proposition 5.4 in [Pop15], reproduced as Proposition 2.10. in [MP21]: a compact complex manifold $X$ is {\it degenerate balanced} if and only if there exists no non-zero $d$-closed $(1, 1)$-current $T\geq 0$ on $X$. In other words, the compact {\it degenerate balanced} manifolds $X$ are characterised by their pseudo-effective cone ${\cal E}_X$ (namely, the set of Bott-Chern cohomology classes of $d$-closed positive $(1, 1)$-currents on $X$) being reduced to the zero class.

This prompts one to ask the following

\begin{Question}\label{Question:small_psef-cone} Let $X$ be a compact complex manifold. Is it true that $X$ is {\bf balanced hyperbolic} if and only if its pseudo-effective cone ${\cal E}_X$ is {\bf small} (in a sense to be determined)?

\end{Question}

In $\S.$\ref{subsection:degree2_DR} and $\S.$\ref{subsection:degree2_BC-A} we provide some evidence for this by first showing that both the {\it balanced} hypothesis on a given Hermitian metric $\omega$ (see Lemma and Definition \ref{Lem-Def:primitive-classes_deg2}) and the {\it Gauduchon} hypothesis (see Lemma and Definition \ref{Lem-Def:primitive-classes_bideg11}) enable one to define a notion of $\omega$-{\it primitive} De Rham cohomology classes of degree $2$ (resp. $\omega$-{\it primitive} Bott-Chern cohomology classes of bidegree $(1,\,1)$). For example, if $\omega$ is {\it balanced}, we set \begin{equation*}H^2_{DR}(X,\,\C)_{\omega\mhyphen prim}:=\ker\bigg(\{\omega_{n-1}\}_{DR}\wedge\cdot\bigg)\subset H^2_{DR}(X,\,\C),\end{equation*} after we have showed that the linear map: \begin{equation*}\{\omega_{n-1}\}_{DR}\wedge\cdot : H^2_{DR}(X,\,\C)\longrightarrow H^{2n}_{DR}(X,\,\C)\simeq\C, \hspace{3ex} \{\alpha\}_{DR}\longmapsto\{\omega_{n-1}\wedge\alpha\}_{DR},\end{equation*} is well defined and depends only on the cohomology class $\{\omega_{n-1}\}_{DR}\in H^{2n-2}_{DR}(X,\,\C)$. We go on to show that a class ${\mathfrak c}\in H^2_{DR}(X,\,\C)$ is $\omega$-primitive if and only if it can be represented by an $\omega$-primitive form (cf. Lemma \ref{Lem:primitive-classes_rep}), a fact that does not seem to hold in the Gauduchon context of $\S.$\ref{subsection:degree2_BC-A}. We then show that, when the balanced metric $\omega$ is {\it not degenerate balanced}, the $\omega$-primitive classes form a complex hyperplane $H^2_{DR}(X,\,\C)_{\omega\mhyphen prim}$ in $H^2_{DR}(X,\,\C)$ that depends only on the balanced class $\{\omega_{n-1}\}_{DR}$. (See Corollary \ref{Cor:primitive_bal_deg-bal}.) Finally, we are able to define a {\it positive side} $H^2_{DR}(X,\,\R)_\omega^{+}$ and a {\it negative side} $H^2_{DR}(X,\,\R)_\omega^{-}$ of the hyperplane $H^2_{DR}(X,\,\R)_{\omega\mhyphen prim}:=H^2_{DR}(X,\,\C)_{\omega\mhyphen prim}\cap H^2_{DR}(X,\,\R)$ in $H^2_{DR}(X,\,\R)$
and get a {\it partition} of $H^2_{DR}(X,\,\R)$: \begin{eqnarray*}H^2_{DR}(X,\,\R) = H^2_{DR}(X,\,\R)_\omega^{+}\cup H^2_{DR}(X,\,\R)_{\omega\mhyphen prim}\cup H^2_{DR}(X,\,\R)_\omega^{-}.\end{eqnarray*}

A similar study of the case where $\omega$ is only a Gauduchon metric in $\S.$\ref{subsection:degree2_BC-A} leads to the characterisation of the pseudo-effective cone as the intersection of the {\it non-negative sides} of all the hyperplanes $H^{1,\,1}_{BC}(X,\,\R)_{\omega\mhyphen prim}$ determined by Aeppli cohomology classes $[\omega_{n-1}]_A$ of Gauduchon metrics $\omega$ on $X$: \begin{equation}\label{eqn:psef-cone_positive-side-char}{\cal E}_X = \bigcap\limits_{[\omega_{n-1}]_A\in{\cal G}_X} H^{1,\,1}_{BC}(X,\,\R)_\omega^{\geq 0},\end{equation} where ${\cal G}_X$ is the {\it Gauduchon cone} of $X$ (introduced in [Pop15] as the set of all such Aeppli cohomology classes, see $\S.$\ref{subsection:deg-bal_background} for a reminder of the terminology).

\vspace{2ex}

In $\S.$\ref{subsection:L1-currents_universal-cover}, we answer a version of Question \ref{Question:small_psef-cone} on the universal covering space of a balanced hyperbolic manifold in the following form. Throughtout the paper, $L^p_{\widetilde\omega}$, $L^p_\omega$, resp. $L^p_g$ will stand for the space of objects that are $L^p$ with respect to the metric $\widetilde\omega$, $\omega$, resp. $g$.

\begin{Prop}\label{Prop:no-L1-currents} Let $(X,\,\omega)$ be a {\bf balanced hyperbolic manifold} and let $\pi:\widetilde{X}\longrightarrow X$ be the universal cover of $X$. There exists no non-zero $d$-closed positive $(1,\,1)$-current $\widetilde{T}\geq 0$ on $\widetilde{X}$ such that $\widetilde{T}$ is $L^1_{\widetilde\omega}$, where $\widetilde\omega:=\pi^\star\omega$ is the lift of $\omega$ to $\widetilde{X}$.

\end{Prop}

This result provides the link with the second part of the paper that we now briefly outline.

\vspace{3ex}

{\bf (II) Case of balanced hyperbolic manifolds}

\vspace{1ex}

The results in the second part of the paper mirror, to some extent, those in the first part. The main difference is that the stage changes from $X$ to its universal covering space $\widetilde{X}$. Specifically, when $X$ is supposed to be balanced hyperbolic, we obtain vanishing theorems for the $L^2$ harmonic cohomology of $\widetilde{X}$.

\vspace{1ex}

{\bf (a)}\, In this setting, our main result in {\bf degree 1} and its dual {\bf degree $2n-1$} is the following

\begin{The}\label{The:bal-complete-exact_no-harmonic-forms_degree1} Let $X$ be a compact complex {\bf balanced hyperbolic} manifold with $\mbox{dim}_\C X=n$. Let $\pi:\widetilde{X}\longrightarrow X$ be the universal cover of $X$ and $\widetilde\omega:=\pi^\star\omega$ the lift to $\widetilde{X}$ of a balanced hyperbolic metric $\omega$ on $X$.

   There are no non-zero $\Delta_{\widetilde\omega}$-harmonic $L^2_{\widetilde\omega}$-forms of pure types and of degrees $1$ and $2n-1$ on $\widetilde{X}$: $${\cal H}^{1,\,0}_{\Delta_{\widetilde\omega}}(\widetilde{X},\,\C) = {\cal H}^{0,\,1}_{\Delta_{\widetilde\omega}}(\widetilde{X},\,\C) = 0  \hspace{3ex} \mbox{and} \hspace{3ex} {\cal H}^{n,\,n-1}_{\Delta_{\widetilde\omega}}(\widetilde{X},\,\C) = {\cal H}^{n-1,\,n}_{\Delta_{\widetilde\omega}}(\widetilde{X},\,\C) = 0,$$ where $\Delta_{\widetilde\omega}:=dd^\star_{\widetilde\omega} + d^\star_{\widetilde\omega}d$ is the $d$-Laplacian induced by the metric $\widetilde\omega$.

\end{The}

The differential operators $d, d^\star_{\widetilde\omega}, \Delta_{\widetilde\omega}$ and all the similar ones are considered as {\it closed} and {\it densely defined} unbounded operators on the spaces $L^2_k(\widetilde{X},\,\C)$ of $L^2_{\widetilde\omega}$-forms of degree $k$ on the {\it complete} complex manifold $(\widetilde{X},\,\widetilde\omega)$. (See reminder of some basic results on complete Riemannian manifolds and unbounded operators in $\S.$\ref{subsection:L1-currents_universal-cover}.)

\vspace{3ex}

{\bf (b)}\, To introduce our results in {\bf degree 2}, we start by reminding the reader of the following facts of [Dem84] (see also [Dem97, VII, $\S.1$]). For any Hermitian metric $\omega$ on a complex manifold $X$ with $\mbox{dim}_\C X=n$, one defines the {\it torsion operator} $\tau=\tau_\omega:=[\Lambda_\omega,\,\partial\omega\wedge\cdot\,]$ of order $0$ and of type $(1,\,0)$ acting on the differential forms of $X$, where $\Lambda_\omega$ is the adjoint of the multiplication operator $\omega\wedge\cdot$ w.r.t. the pointwise inner product $\langle\,\,,\,\,\rangle_\omega$ defined by $\omega$. The K\"ahler commutation relations generalise to the arbitrary Hermitian context as \begin{eqnarray}\label{eqn:Hermitian-commutation}i[\Lambda_\omega,\,\bar\partial] = \partial^\star + \tau^\star\end{eqnarray} and the three other relations obtained from this one by conjugation and/or adjunction. (See [Dem97, VII, $\S.1$, Theorem 1.1.].) Moreover, considering the torsion-twisted Laplacians $$\Delta_\tau:=[d+\tau,\,d^\star + \tau^\star]  \hspace{3ex} \mbox{and} \hspace{3ex} \Delta'_\tau:=[\partial+\tau,\,\partial^\star + \tau^\star],$$ the following formula holds (see [Dem97, VII, $\S.1$, Proposition 1.16.]): \begin{eqnarray}\label{eqn:Delta_tau_sum-formula}\Delta_\tau = \Delta'_\tau + \Delta''.\end{eqnarray} When the metric $\omega$ is K\"ahler, one has $\tau=0$ and one recovers the classical formula $\Delta = \Delta' + \Delta''$. However, we will deal with a more general, possibly non-K\"ahler, case.

\vspace{3ex}

In the context of balanced hyperbolic manifolds, our main result in {\bf degree 2} is the following

\begin{The}\label{The:Delta_tau_harmonic_bal-complete-exact} Let $X$ be a compact complex {\bf balanced hyperbolic} manifold with $\mbox{dim}_\C X=n$. Let $\pi:\widetilde{X}\longrightarrow X$ be the universal cover of $X$ and $\widetilde\omega:=\pi^\star\omega$ the lift to $\widetilde{X}$ of a balanced hyperbolic metric $\omega$ on $X$.

   There are no non-zero semi-positive $\Delta_{\widetilde\tau}$-harmonic $L^2_{\widetilde\omega}$-forms of pure type $(1,\,1)$ on $\widetilde{X}$: $$\bigg\{\alpha^{1,\,1}\in{\cal H}^{1,\,1}_{\Delta_{\widetilde\tau}}(\widetilde{X},\,\C)\,\mid\,\alpha^{1,\,1}\geq 0\bigg\} = \{0\},$$ where $\widetilde\tau=\widetilde\tau_{\widetilde\omega}:=[\Lambda_{\widetilde\omega},\,\partial\widetilde\omega\wedge\cdot]$.

\end{The}

\vspace{2ex}

As a piece of notation that will be used throughout the text, whenever $u$ is a $k$-form and $(p,\,q)$ is a bidegree with $p+q=k$, $u^{p,\,q}$ will stand for the component of $u$ of bidegree $(p,\,q)$. 

\vspace{2ex}

\noindent {\bf Acknowledgments.} This work is part of the first-named author's PhD thesis under the supervision of the second-named author. The former wishes to express his gratitude to the latter for his constant guidance while this work was carried out, as well as to his Tunisian supervisor, Fathi Haggui, for constant support. Both authors are very grateful to the referee for their careful reading of the manuscript and their useful remarks and suggestions.

\section{Properties of degenerate balanced manifolds}\label{section:deg-bal_prop}

In this section, we investigate the effect of the balanced condition on the cohomology of degrees $1$ and $2$, while pointing out the peculiarities of the degenerate balanced case.

\subsection{Background}\label{subsection:deg-bal_background}

Given a complex manifold $X$ with $\mbox{dim}_\C X=n\geq 2$ and a Hermitian metric $\omega$ (identified with its underlying $C^\infty$ positive definite $(1,\,1)$-form $\omega$) on $X$, we will put $\omega_r:=\omega^r/r!$ for $r=1,\dots, n$. Moreover, we denote by $C^{r,\,s}(X) = C^{r,\,s}(X,\,\C)$ the space of smooth $\C$-valued $(r,\,s)$-forms on $X$ for $r,s=1,\dots , n$. If $X$ is compact, recall the classical definitions of the {\it Bott-Chern} and {\it Aeppli} cohomology groups of $X$ of any bidegree $(p,\,q)$: \begin{eqnarray*}H^{p,\,q}_{BC}(X,\,\C) = \frac{\ker(\partial:C^{p,\,q}(X)\to C^{p+1,\,q}(X))\cap\ker(\bar\partial:C^{p,\,q}(X)\to C^{p,\,q+1}(X))}{\mbox{Im}\,(\partial\bar\partial:C^{p-1,\,q-1}(X)\to C^{p,\,q}(X))}\end{eqnarray*}
\begin{eqnarray*}  H^{p,\,q}_A(X,\,\C) = \frac{\ker(\partial\bar\partial:C^{p,\,q}(X)\to C^{p+1,\,q+1}(X))}{\mbox{Im}\,(\partial:C^{p-1,\,q}(X)\to C^{p,\,q}(X)) + \mbox{Im}\,(\bar\partial:C^{p,\,q-1}(X)\to C^{p,\,q}(X))}.\end{eqnarray*}

We will use the Serre-type duality (see e.g. [Sch07]): \begin{equation}\label{eqn:BC-A_duality}H^{1,\,1}_{BC}(X,\,\C)\times H^{n-1,\,n-1}_A(X,\,\C)\longrightarrow\C, \hspace{3ex} (\{u\}_{BC},\,\{v\}_A)\mapsto\{u\}_{BC}.\{v\}_A:=\int\limits_Xu\wedge v,\end{equation} as well as the {\it pseudo-effective cone} of $X$ introduced in [Dem92, Definition 1.3] as the set $${\cal E}_X:=\bigg\{[T]_{BC}\in H^{1,\,1}_{BC}(X,\,\R)\,\slash\, T\geq 0\hspace{1ex} d\mbox{-closed} \hspace{1ex} (1,\,1)\mbox{-current} \hspace{1ex} \mbox{on} \hspace{1ex} X\bigg\}.$$

Recall that a Hermitian metric $\omega$ on $X$ is said to be a {\it Gauduchon metric} (cf. [Gau77]) if $\partial\bar\partial\omega^{n-1}=0$. For any such metric $\omega$, $\omega^{n-1}$ defines an Aeppli cohomology class and the set of all these cohomology classes is called the {\it Gauduchon cone} of $X$ (cf. [Pop15]): $${\cal G}_X:=\bigg\{\{\omega^{n-1}\}_A\in H^{n-1,\,n-1}_A(X,\,\R)\,\mid\,\omega\hspace{1ex}\mbox{is a Gauduchon metric}\hspace{1ex} \mbox{on}\hspace{1ex}X\bigg\}\subset H^{n-1,\,n-1}_A(X,\,\R).$$

The main link between the cones ${\cal G}_X$ and ${\cal E}_X$ on a compact $n$-dimensional $X$ is provided by the following reformulation observed in [Pop16] of a result of Lamari's from [Lam99]. The pseudo-effective cone ${\cal E}_X\subset H^{1,\,1}_{BC}(X,\,\R)$ and the closure of the Gauduchon cone $\overline{\cal G}_X\subset H^{n-1,\,n-1}_A(X,\,\R)$ are {\bf dual} to each other under the duality (\ref{eqn:BC-A_duality}). This means that the following two statements hold.

\vspace{1ex}

$(1)$\, Given any class $\mathfrak{c}^{1,\,1}_{BC}\in H^{1,\,1}_{BC}(X,\,\R)$, the following equivalence holds: $$\mathfrak{c}^{1,\,1}_{BC}\in{\cal E}_X\iff\mathfrak{c}^{1,\,1}_{BC}.\mathfrak{c}^{n-1,\,n-1}_A\geq 0 \hspace{2ex} \mbox{for every class} \hspace{1ex} \mathfrak{c}^{n-1,\,n-1}_A\in{\cal G}_X.$$

\vspace{1ex}

$(2)$\, Given any class $\mathfrak{c}^{n-1,\,n-1}_A\in H^{n-1,\,n-1}_A(X,\,\R)$, the following equivalence holds: $$\mathfrak{c}^{n-1,\,n-1}_A\in\overline{\cal G}_X\iff\mathfrak{c}^{1,\,1}_{BC}.\mathfrak{c}^{n-1,\,n-1}_A\geq 0 \hspace{2ex} \mbox{for every class} \hspace{1ex} \mathfrak{c}^{1,\,1}_{BC}\in{\cal E}_X.$$ 

\vspace{2ex}

Finally, recall that a compact complex manifold $X$ is said to be a $\partial\bar\partial$-{\bf manifold} (see [DGMS75] for the notion, [Pop14] for the name) if for any $d$-closed {\it pure-type} form $u$ on $X$, the following exactness properties are equivalent:

\vspace{1ex}

\hspace{10ex} $u$ is $d$-exact $\Longleftrightarrow$ $u$ is $\partial$-exact $\Longleftrightarrow$ $u$ is $\bar\partial$-exact $\Longleftrightarrow$ $u$ is $\partial\bar\partial$-exact.

\vspace{2ex}
  
On a complex manifold $X$ with $\mbox{dim}_\C X=n$, we will often use the following standard formula (cf. e.g. [Voi02, Proposition 6.29, p. 150]) for the Hodge star operator $\star = \star_\omega$ of any Hermitian metric $\omega$ applied to $\omega$-{\it primitive} forms $v$ of arbitrary bidegree $(p, \, q)$: \begin{eqnarray}\label{eqn:prim-form-star-formula-gen}\star\, v = (-1)^{k(k+1)/2}\, i^{p-q}\,\omega_{n-p-q}\wedge v, \hspace{2ex} \mbox{where}\,\, k:=p+q.\end{eqnarray} Recall that, for any $k=0,1, \dots , n$, a $k$-form $v$ is said to be $(\omega)$-{\it primitive} if $\omega_{n-k+1}\wedge v=0$ and that this condition is equivalent to $\Lambda_\omega v=0$, where $\Lambda_\omega$ is the adjoint of the operator $\omega\wedge\cdot$ (of multiplication by $\omega$) w.r.t. the pointwise inner product $\langle\,\,,\,\,\rangle_\omega$ defined by $\omega$.

We will also often deal with $C^\infty$ $(1,\,1)$-forms $\alpha$. If $\alpha = \alpha_{prim} + f\omega$ is the {\it Lefschetz decomposition}, where $\alpha_{prim}$ is {\it primitive} and $f$ is a smooth function on $X$, we get $\Lambda_\omega\alpha = nf$, hence \begin{eqnarray}\label{eqn:Lefschetz_decomp_11}\alpha = \alpha_{prim} + \frac{1}{n}\,(\Lambda_\omega\alpha)\,\omega.\end{eqnarray} We will often write $(1,\,1)$-forms in this form.

On the other hand, we will often indicate the metric with respect to which certain operators are calculated. For example, $d^\star_\omega$ and $\Delta_\omega:=dd^\star_\omega + d^\star_\omega d$ are the adjoint of $d$, resp. the $d$-Laplacian, induced by the metric $\omega$.

\subsection{Case of degree $1$}\label{subsection:degree1} The starting point is the following

\begin{Lem}\label{Lem:1-forms_Delta-harm} Let $\omega$ be a Hermitian metric on a complex manifold $X$ with $\mbox{dim}_\C X=n$. Fix a form $u = u^{1,\,0} + u^{0,\,1}\in C^\infty_1(X,\,\C)$.

  \vspace{1ex}

  (i)\, The following formula holds: \begin{eqnarray}\label{eqn:1-forms_d-star_omega_n-1}\nonumber d^\star(\omega_{n-1}\wedge u) & = & i(\partial u^{1,\,0} - \bar\partial u^{0,\,1})\wedge\omega_{n-2} + i\bigg((\partial u^{0,\,1})_{prim} - (\bar\partial u^{1,\,0})_{prim}\bigg)\wedge\omega_{n-2}\\
    & + & \frac{i}{n}\,\bigg(\Lambda_\omega(\bar\partial u^{1,\,0}) - \Lambda_\omega(\partial u^{0,\,1})\bigg)\,\omega_{n-1},\end{eqnarray} where $d^\star = d^\star_\omega$ is the formal adjoint of $d$ w.r.t. the $L^2_\omega$ inner product, while the subscript {\rm prim} indicates the $\omega$-primitive part in the Lefschetz decomposition of the form to which it is applied.

  In particular, if $du^{1,\,0} = 0$ and $du^{0,\,1} = 0$, we get \begin{eqnarray*}d^\star(\omega_{n-1}\wedge u) = 0.\end{eqnarray*}

 (ii)\, If $\omega$ is {\bf balanced} and $du^{1,\,0} = 0$ and $du^{0,\,1} = 0$, then \begin{eqnarray*}\Delta(\omega_{n-1}\wedge u) = 0,\end{eqnarray*} where $\Delta = \Delta_\omega = dd^\star + d^\star d$ is the $d$-Laplacian induced by $\omega$.

 (iii)\, If $X$ is compact, $\omega$ is {\bf degenerate balanced} and $du^{1,\,0} = du^{0,\,1} = 0$, then $u=0$.
  
\end{Lem}

\noindent {\it Proof.} (i)\, All $1$-forms are primitive, so from the standard formula (\ref{eqn:prim-form-star-formula-gen}) we get: $\star\, u^{1,\,0} = -i\omega_{n-1}\wedge u^{1,\,0}$, hence $\star\,(\omega_{n-1}\wedge u^{1,\,0}) = -iu^{1,\,0}$. Meanwhile, $d^\star = -\star d\star$, so applying $-\star d$ to the previous identity and writing the $(1,\,1)$-form $\bar\partial u^{1,\,0}$ in the form (\ref{eqn:Lefschetz_decomp_11}), we get the first line below: \begin{eqnarray*}d^\star(\omega_{n-1}\wedge u^{1,\,0}) & = & i\,\star\partial u^{1,\,0} + i\,\star(\bar\partial u^{1,\,0})_{prim} + \frac{i}{n}\,(\Lambda_\omega\bar\partial u^{1,\,0})\,\star\omega \\
  & = & i\partial u^{1,\,0}\wedge\omega_{n-2} - i(\bar\partial u^{1,\,0})_{prim}\wedge\omega_{n-2} + \frac{i}{n}\,(\Lambda_\omega\bar\partial u^{1,\,0})\,\omega_{n-1},\end{eqnarray*} where the second line follows from the standard formula (\ref{eqn:prim-form-star-formula-gen}).

Running the analogous computations for $u^{0,\,1}$ or taking conjugates, we get:

\begin{eqnarray*}d^\star(\omega_{n-1}\wedge u^{0,\,1}) = -i\bar\partial u^{0,\,1}\wedge\omega_{n-2} + i(\partial u^{0,\,1})_{prim}\wedge\omega_{n-2} - \frac{i}{n}\,(\Lambda_\omega\partial u^{0,\,1})\,\omega_{n-1}.\end{eqnarray*}

Formula (\ref{eqn:1-forms_d-star_omega_n-1}) follows by adding up the above expressions for $d^\star(\omega_{n-1}\wedge u^{1,\,0})$ and $d^\star(\omega_{n-1}\wedge u^{0,\,1})$.

\vspace{1ex}

(ii)\, If $\omega$ is balanced, we get $d(\omega_{n-1}\wedge u) = \omega_{n-1}\wedge du =0$ since $du=0$ under the assumptions. Since we also have $d^\star(\omega_{n-1}\wedge u)=0$ by (i), the contention follows.

\vspace{1ex}

(iii)\, If $\omega$ is degenerate balanced, there exists a smooth $(2n-3)$-form $\Gamma$ such that $\omega^{n-1} = d\Gamma$. Hence, $\omega^{n-1}\wedge u = d(\Gamma\wedge u)\in\mbox{Im}\,d$ because we also have $du=0$ by our assumptions. However, $\omega^{n-1}\wedge u\in\ker\Delta$ by (ii) and $\ker\Delta\perp\mbox{Im}\,d$ by the compactness assumption on $X$. Thus, the form $\omega^{n-1}\wedge u\in\ker\Delta\cap\mbox{Im}\,d = \{0\}$ must vanish.

On the other hand, the pointwise map $\omega_{n-1}\wedge\cdot : \Lambda^1T^\star X\longrightarrow\Lambda^{2n-1}T^\star X$ is bijective, so we get $u=0$ from $\omega^{n-1}\wedge u=0$.  \hfill $\Box$

\vspace{3ex}

We now use Lemma \ref{Lem:1-forms_Delta-harm} to infer its consequences announced in the introduction.

\vspace{2ex}

\noindent $\bullet$ {\bf Proof of (i) of Proposition \ref{Prop:BC_10-01_b1_deg-bal}.} This follows at once from (iii) of Lemma \ref{Lem:1-forms_Delta-harm}.  \hfill $\Box$

\vspace{2ex}

Another consequence of Lemma \ref{Lem:1-forms_Delta-harm} is that the balanced condition, combined with the mild $\partial\bar\partial$-type condition in (ii) of Theorem \ref{The:Hard-Lefschetz_bal-ddbar}, enables one to get a {\it Hard Lefschetz Isomorphism} between the De Rham cohomology spaces of degrees $1$ and $2n-1$.

\vspace{2ex}

\noindent $\bullet$ {\bf Proof of Theorem \ref{The:Hard-Lefschetz_bal-ddbar}.} (i)\, Lemma \ref{Lem:1-forms_Delta-harm} is not needed here. Let $u$ be a smooth $1$-form. Since $d\omega_{n-1}=0$, $d(\omega_{n-1}\wedge u)=0$ whenever $du=0$, while $\omega_{n-1}\wedge u=d(f\omega_{n-1})$ whenever $u=df$ for some smooth function $f$ on $X$. This proves the well-definedness of the map (\ref{eqn:Hard-Lefschetz_map}).

Similarly, if $\omega_{n-1} = \gamma_{n-1} + d\Gamma$ for some smooth $(2n-2)$-form $\gamma_{n-1}$ and some smooth $(2n-3)$-form $\Gamma$, then $\omega_{n-1}\wedge u = \gamma_{n-1}\wedge u + d(\Gamma\wedge u)$ for every $d$-closed $1$-form $u$. Hence, $\{\omega_{n-1}\wedge u\}_{DR} = \{\gamma_{n-1}\wedge u\}_{DR}$ whenever $\{\omega_{n-1}\}_{DR} = \{\gamma_{n-1}\}_{DR}$, so the map (\ref{eqn:Hard-Lefschetz_map}) depends only on $\{\omega_{n-1}\}_{DR}$.

\vspace{1ex}

(ii)\, Since $H^1_{DR}(X,\,\C)$ and $H^{2n-1}_{DR}(X,\,\C)$ have equal dimensions, by Poincar\'e duality, it suffices to prove that the map (\ref{eqn:Hard-Lefschetz_map}) is injective.

Let $u$ be an arbitrary smooth $d$-closed $1$-form on $X$. We start by showing that there exists a smooth function $f:X\to\C$ such that $\partial u^{0,\,1} = \partial\bar\partial f$ on $X$. To see this, notice that the property $du=0$ translates to the following three relations holding: \begin{equation}\label{eqn:Hard-Lefschetz_bal-ddbar_proof1}(a)\,\,\partial u^{1,\,0} = 0; \hspace{3ex} (b)\,\, \partial u^{0,\,1} + \bar\partial u^{1,\,0} = 0; \hspace{3ex} (c)\,\, \bar\partial u^{0,\,1} = 0.\end{equation} Thus, the $(1,\,1)$-form $\partial u^{0,\,1}$ is $d$-closed (since it is $\bar\partial$-closed by (c) of (\ref{eqn:Hard-Lefschetz_bal-ddbar_proof1})) and $\partial$-exact. Thanks to hypothesis (\ref{eqn:Hard-Lefschetz_bal-ddbar}), we infer that $\partial u^{0,\,1}$ is $\partial\bar\partial$-exact. Thus, there exists a smooth function $f$ as stated.

Using (b) of (\ref{eqn:Hard-Lefschetz_bal-ddbar_proof1}), we further infer that $\bar\partial u^{1,\,0} = -\partial u^{0,\,1} = -\partial\bar\partial f$, so $\bar\partial(u^{1,\,0} - \partial f) = 0$. From the identities $\partial(u^{0,\,1} - \bar\partial f) = 0$ and $\bar\partial(u^{1,\,0} - \partial f) = 0$ and from (a) and (c) of (\ref{eqn:Hard-Lefschetz_bal-ddbar_proof1}), we get: $$d(u^{1,\,0} - \partial f) = 0 \hspace{3ex} \mbox{and} \hspace{3ex} d(u^{0,\,1} - \bar\partial f) = 0.$$ This means that $$(u - df)^{1,\,0}\in\ker d \hspace{3ex} \mbox{and} \hspace{3ex} (u - df)^{0,\,1}\in\ker d.$$ From this and from (i) of Lemma \ref{Lem:1-forms_Delta-harm}, we deduce that \begin{equation}\label{eqn:Hard-Lefschetz_bal-ddbar_proof2}\omega_{n-1}\wedge(u - df)\in\ker d^\star.\end{equation}

On the other hand, if $\omega$ is balanced and if $\{\omega_{n-1}\wedge u\}_{DR}=0\in H^{2n-1}_{DR}(X,\,\C)$ (i.e. $\omega_{n-1}\wedge u\in\mbox{Im}\,d$), then \begin{equation}\label{eqn:Hard-Lefschetz_bal-ddbar_proof3}\omega_{n-1}\wedge(u - df)\in\mbox{Im}\,d.\end{equation}

From (\ref{eqn:Hard-Lefschetz_bal-ddbar_proof2}), (\ref{eqn:Hard-Lefschetz_bal-ddbar_proof3}) and $\ker d^\star\perp\mbox{Im}\,d$, we infer that $\omega_{n-1}\wedge(u - df)=0$. Since $u - df$ is a smooth $1$-form on $X$ and the pointwise-defined linear map: $$\omega_{n-1}\wedge\cdot : C^\infty_1(X,\,C)\longrightarrow C^\infty_{2n-1}(X,\,C), \hspace{3ex} \alpha\mapsto\omega_{n-1}\wedge\alpha,$$ is bijective, we finally get $u - df=0$, so $\{u\}_{DR}=0\in H^1_{DR}(X,\,\C)$.

This proves the injectivity of the map (\ref{eqn:Hard-Lefschetz_map}) whenever $\omega$ is balanced and $X$ satisfies hypothesis (\ref{eqn:Hard-Lefschetz_bal-ddbar}).  \hfill $\Box$

\vspace{2ex}

In the degenerate balanced case, we get the vanishing of the first Betti number of the manifold.

\vspace{2ex}

\noindent $\bullet$ {\bf Proof of (ii) of Proposition \ref{Prop:BC_10-01_b1_deg-bal}.} If $\omega$ is degenerate balanced, the map (\ref{eqn:Hard-Lefschetz_map}) vanishes identically. Meanwhile, by Theorem \ref{The:Hard-Lefschetz_bal-ddbar}, the map (\ref{eqn:Hard-Lefschetz_map}) is an isomorphism. We get $H^1_{DR}(X,\,\C)=0$. \hfill $\Box$

\subsection{Case of degree $2$: De Rham cohomology}\label{subsection:degree2_DR} The balanced property of a metric enables one to define a notion of primitivity for $2$-forms.

\begin{Lem-Def}\label{Lem-Def:primitive-classes_deg2} Let $\omega$ be a {\bf balanced} metric on a compact complex manifold $X$ with $\mbox{dim}_\C X=n$. The linear map: \begin{equation}\label{eqn:primitive-classes_deg2}\{\omega_{n-1}\}_{DR}\wedge\cdot : H^2_{DR}(X,\,\C)\longrightarrow H^{2n}_{DR}(X,\,\C)\simeq\C, \hspace{3ex} \{\alpha\}_{DR}\longmapsto\{\omega_{n-1}\wedge\alpha\}_{DR},\end{equation} is {\bf well defined} and depends only on the cohomology class $\{\omega_{n-1}\}_{DR}\in H^{2n-2}_{DR}(X,\,\C)$. We set: \begin{equation*}H^2_{DR}(X,\,\C)_{\omega\mhyphen prim}:=\ker\bigg(\{\omega_{n-1}\}_{DR}\wedge\cdot\bigg)\subset H^2_{DR}(X,\,\C)\end{equation*} and we call its elements {\bf ($\omega$-)primitive} De Rham $2$-classes.

\end{Lem-Def}

\noindent {\it Proof.} Since $d\omega_{n-1}=0$, for every $d$-closed (resp. $d$-exact) $2$-form $\alpha$, $\omega_{n-1}\wedge\alpha$ is $d$-closed (resp. $d$-exact). This proves the well-definedness of the map.

Meanwhile, if $\Omega\in C^\infty_{n-1,\,n-1}(X,\,\C)$ is such that $\Omega = \omega_{n-1} + d\Gamma$ for some smooth $(2n-3)$-form $\Gamma$, then, for every $d$-closed $2$-form $\alpha$, $\Omega\wedge\alpha = \omega_{n-1}\wedge\alpha + d(\Gamma
\wedge\alpha)$. Hence, $\{\Omega\wedge\alpha\}_{DR} = \{\omega_{n-1}\wedge\alpha\}_{DR}$ whenever $\{\Omega\}_{DR} = \{\omega_{n-1}\}_{DR}$. This proves the independence of the map $\{\omega_{n-1}\}_{DR}\wedge\cdot$ of the choice of representative of the class $\{\omega_{n-1}\}_{DR}$. \hfill $\Box$

\vspace{2ex}

We now observe a link between {\it primitive} $2$-classes and {\it primitive} $2$-forms.

\begin{Lem}\label{Lem:primitive-classes_rep} Let $\omega$ be a {\bf balanced} metric on a compact complex manifold $X$ with $\mbox{dim}_\C X=n$. For any class ${\mathfrak c}\in H^2_{DR}(X,\,\C)$, the following equivalence holds: $${\mathfrak c} \hspace{1ex} \mbox{is {\bf $\omega$-primitive}} \iff \exists\,\alpha\in{\mathfrak c} \hspace{2ex} \mbox{such that} \hspace{2ex}\alpha \hspace{1ex} \mbox{is {\bf $\omega$-primitive}}.$$

\end{Lem}

\noindent {\it Proof.} ``$\Longleftarrow$'' Suppose $\alpha\in C^\infty_2(X,\,\C)$ such that $d\alpha=0$, $\alpha\in{\mathfrak c}$ and $\alpha$ is $\omega$-primitive. Then, $\omega_{n-1}\wedge\alpha=0$, hence $\{\omega_{n-1}\wedge\alpha\}_{DR}=0$. This means that the class ${\mathfrak c}=\{\alpha\}_{DR}$ is $\omega$-primitive.

\vspace{1ex}

``$\implies$'' Suppose the class ${\mathfrak c}$ is $\omega$-primitive. Pick an arbitrary representative $\beta\in{\mathfrak c}$. The $\omega$-primitivity of ${\mathfrak c}=\{\beta\}_{DR}$ translates to $\{\omega_{n-1}\wedge\beta\}_{DR} = 0\in H^{2n}_{DR}(X,\,\C)$. This, in turn, is equivalent to the existence of a form $\Gamma\in C^\infty_{2n-1}(X,\,\C)$ such that $\omega_{n-1}\wedge\beta = d\Gamma$.

Meanwhile, we know from the general theory that the map $$\omega_{n-1}\wedge\cdot : C^\infty_1(X,\,\C)\longrightarrow C^\infty_{2n-1}(X,\,\C)$$ is an isomorphism. Hence, there exists a unique $u\in C^\infty_1(X,\,\C)$ such that $\Gamma = \omega_{n-1}\wedge u$. We get: \begin{eqnarray*}\omega_{n-1}\wedge\beta = d\Gamma = \omega_{n-1}\wedge du,\end{eqnarray*} where the last identity follows from the balanced property of $\omega$. Consequently, $$\omega_{n-1}\wedge(\beta -du) = 0,$$ proving that $\alpha:=\beta -du$ is a {\it primitive} representative of the class ${\mathfrak c}=\{\beta\}_{DR}$.  \hfill $\Box$

\vspace{2ex}

Finally, we can characterise the {\it degenerate balanced} property of a given {\it balanced} metric in terms of primitivity for $2$-classes.

\begin{Lem}\label{Lem:primitive_deg-bal_char} Let $\omega$ be a {\bf balanced} metric on a compact complex manifold $X$ with $\mbox{dim}_\C X=n$. The following equivalence holds: $$H^2_{DR}(X,\,\C)_{\omega\mhyphen prim}= H^2_{DR}(X,\,\C) \iff \omega \hspace{1ex} \mbox{is {\bf degenerate balanced}}.$$

\end{Lem}

\noindent {\it Proof.} ``$\Longleftarrow$'' Suppose that $\omega$ is degenerate balanced. Then $\omega_{n-1}$ is $d$-exact, hence $\omega_{n-1}\wedge\alpha$ is $d$-exact (or equivalently $\{\omega_{n-1}\wedge\alpha\}_{DR} = 0\in H^{2n}_{DR}(X,\,\C)$) for every $d$-closed $2$-form $\alpha$. This means that the map $\{\omega_{n-1}\}_{DR}\wedge\cdot : H^2_{DR}(X,\,\C)\longrightarrow H^{2n}_{DR}(X,\,\C)$ vanishes identically, so $H^2_{DR}(X,\,\C)_{\omega\mhyphen prim}= H^2_{DR}(X,\,\C)$.

\vspace{1ex}

``$\implies$'' Suppose that $H^2_{DR}(X,\,\C)_{\omega\mhyphen prim}= H^2_{DR}(X,\,\C)$. This translates to \begin{equation}\label{eqn:primitive_deg-bal_char_proof_1}\omega_{n-1}\wedge\alpha\in\mbox{Im}\,d, \hspace{6ex}\forall\,\alpha\in C^\infty_2(X,\,\C)\cap\ker d.\end{equation}

Since both $\omega_{n-1}$ and $\alpha$ are $d$-closed, they both have unique $L^2_\omega$-orthogonal decompositions: $$\omega_{n-1} = (\omega_{n-1})_h + d\Gamma \hspace{5ex} \mbox{and} \hspace{5ex} \alpha = \alpha_h + du,$$ where $(\omega_{n-1})_h$ and $\alpha_h$ are $\Delta_\omega$-harmonic, while $\Gamma$ and $u$ are smooth forms of respective degrees $2n-3$ and $1$. We get: $$\omega_{n-1}\wedge\alpha = (\omega_{n-1})_h\wedge\alpha_h + d\bigg((\omega_{n-1})_h\wedge u + \Gamma\wedge\alpha_h + \Gamma\wedge du\bigg), \hspace{6ex}\forall\,\alpha\in C^\infty_2(X,\,\C)\cap\ker d.$$ Together with (\ref{eqn:primitive_deg-bal_char_proof_1}), this implies that \begin{equation}\label{eqn:primitive_deg-bal_char_proof_2}(\omega_{n-1})_h\wedge\alpha_h\in\mbox{Im}\,d, \hspace{6ex}\forall\,\alpha_h\in\ker\Delta_\omega\cap C^\infty_2(X,\,\C).\end{equation}

Meanwhile, since $(\omega_{n-1})_h$ is $\Delta_\omega$-harmonic (and real), $\star_\omega(\omega_{n-1})_h$ is again $\Delta_\omega$-harmonic (and real). Hence, \begin{equation*}\mbox{Im}\,d\ni(\omega_{n-1})_h\wedge\star_\omega(\omega_{n-1})_h = |(\omega_{n-1})_h|^2_\omega\,dV_\omega\geq 0,\end{equation*} where the first relation follows from (\ref{eqn:primitive_deg-bal_char_proof_2}) by choosing $\alpha_h=\star_\omega(\omega_{n-1})_h$. Consequently, from Stokes's Theorem we get: $$\int\limits_X |(\omega_{n-1})_h|^2_\omega\,dV_\omega = 0,$$ hence $(\omega_{n-1})_h=0$. This implies that $\omega_{n-1}$ is $d$-exact, which means that $\omega$ is degenerate balanced.  \hfill $\Box$

\begin{Cor}\label{Cor:primitive_bal_deg-bal} Let $\omega$ be a {\bf balanced} metric on a compact complex manifold $X$ with $\mbox{dim}_\C X=n$. The following dichotomy holds:

  \vspace{1ex}

  (a)\, if $\omega$ is {\bf not degenerate balanced}, $H^2_{DR}(X,\,\C)_{\omega\mhyphen prim}$ is a complex hyperplane in $H^2_{DR}(X,\,\C)$ depending only on the balanced class $\{\omega_{n-1}\}_{DR}$;

  \vspace{1ex}

  (b)\, if $\omega$ is {\bf degenerate balanced}, $H^2_{DR}(X,\,\C)_{\omega\mhyphen prim}= H^2_{DR}(X,\,\C)$.
  
\end{Cor}

\noindent {\it Proof.} This follows immediately from Lemma and Definition \ref{Lem-Def:primitive-classes_deg2}, from Lemma \ref{Lem:primitive_deg-bal_char} and from $H^{2n}_{DR}(X,\,\C)\simeq\C$.  \hfill $\Box$

\vspace{2ex}

We shall now get a {\it Lefschetz-type decomposition} of $H^2_{DR}(X,\,\C)$, induced by an arbitrary balanced metric $\omega$, with $H^2_{DR}(X,\,\C)_{\omega\mhyphen prim}$ as a direct factor. Recall that the balanced condition $d\omega^{n-1}=0$ is equivalent to $d_\omega^\star\omega=0$.

Thanks to the orthogonal $3$-space decompositions: $$C^\infty_k(X,\,\C) = \ker\Delta_\omega\oplus\mbox{Im}\,d\oplus\mbox{Im}\,d_\omega^\star, \hspace{6ex} k\in\{0,\dots , 2n\},$$ where $\ker\Delta_\omega\oplus\mbox{Im}\,d = \ker d$ and $\ker\Delta_\omega\oplus\mbox{Im}\,d_\omega^\star = \ker d_\omega^\star$, applied with $k=2$ and $k=2n-2$, we get unique decompositions of $\omega$, resp. $\omega_{n-1}$: \begin{equation}\label{eqn:omega_power_3-space-decomp}\ker d_\omega^\star\ni\omega = \omega_h + d_\omega^\star\eta_\omega \hspace{3ex} \mbox{and} \hspace{3ex} \ker d\ni\omega_{n-1} = (\omega_{n-1})_h + d\Gamma_\omega,\end{equation} where $\omega_h\in\ker\Delta_\omega$ as a $2$-form, $(\omega_{n-1})_h\in\ker\Delta_\omega$ as a $(2n-2)$-form, while $\eta_\omega$ and $\Gamma_\omega$ are smooth forms of respective degrees $3$ and $2n-3$. Since $\omega$ and $\omega_{n-1}$ are real, so are their harmonic components $\omega_h$ and $(\omega_{n-1})_h$. 

Moreover, it is well known that $\star_\omega\omega = \omega_{n-1}$ and that the Hodge star operator $\star_\omega$ maps $d$-exact forms to $d^\star_\omega$-exact forms and vice-versa. Hence, we get: \begin{equation}\label{eqn:omega_power_3-star-interchange}\star_\omega\omega_h = (\omega_{n-1})_h \hspace{3ex} \mbox{and} \hspace{3ex} \star_\omega(d_\omega^\star\eta_\omega) = d\Gamma_\omega.\end{equation} 

Thus, $\omega_h$ is uniquely determined by $\omega$ and is $d$-closed (because it is even $\Delta_\omega$-harmonic). Therefore, it represents a class in $H^2_{DR}(X,\,\R)$.

\begin{Def}\label{Def:H2-class_bal-metric} For any balanced metric $\omega$ on a compact complex manifold $X$, the De Rham cohomlogy class $\{\omega_h\}_{DR}\in  H^2_{DR}(X,\,\R)$ is called the {\bf cohomology class} of $\omega$.

\end{Def}  

Of course, if $\omega$ is K\"ahler, $\omega_h = \omega$, so $\{\omega_h\}_{DR}$ is the usual K\"ahler class $\{\omega\}_{DR}$.

\begin{Lem}\label{Lem:orthogonality_prim-class-rep_balanced} Suppose there exists a balanced metric $\omega$ on a compact complex manifold $X$. Then, for every $\alpha\in C^\infty_2(X,\,\C)$ such that $d\alpha=0$ and $\{\alpha\}_{DR}\in H^2_{DR}(X,\,\C)_{\omega\mhyphen prim}$, we have: $$\langle\langle\omega_h,\,\alpha\rangle\rangle_\omega = 0,$$ where $\langle\langle\,\,,\,\,\rangle\rangle_\omega$ is the $L^2$ inner product induced by $\omega$.

\end{Lem}

\noindent {\it Proof.} Since $\{\alpha\}_{DR}\in H^2_{DR}(X,\,\C)_{\omega\mhyphen prim}$, there exists a form $\Omega\in C^\infty_{2n-1}$ such that $\omega_{n-1}\wedge\alpha = du$. We get: \begin{eqnarray*}\langle\langle\alpha,\,\omega_h\rangle\rangle_\omega & = & \int\limits_X\alpha\wedge\star_\omega\omega_h \stackrel{(a)}{=} \int\limits_X\alpha\wedge(\omega_{n-1})_h \stackrel{(b)}{=} \int\limits_X\alpha\wedge(\omega_{n-1} - d\Gamma_\omega) \\
  & = & \int\limits_X\alpha\wedge\omega_{n-1} = \int\limits_Xdu = 0,\end{eqnarray*} where Stokes implies two of the last three equalities (note that $\alpha\wedge d\Gamma_\omega = d(\alpha\wedge\Gamma_\omega)$), while (a) follows from (\ref{eqn:omega_power_3-star-interchange}) and (b) follows from (\ref{eqn:omega_power_3-space-decomp}). \hfill $\Box$

\begin{Conc}\label{Conc:Lefschetz-type-decomp_bal} Let $X$ be a compact complex manifold with $\mbox{dim}_\C X=n$. Suppose there exists a non-degenerate {\bf balanced metric} $\omega$ on $X$. Then, the De Rham cohomology space of degree $2$ has a {\bf Lefschetz-type} $L^2_\omega$-orthogonal {\bf decomposition}: \begin{equation}\label{eqn:Lefschetz-type-decomp_bal}H^2_{DR}(X,\,\C) = H^2_{DR}(X,\,\C)_{\omega\mhyphen prim}\oplus\C\cdot\{\omega_h\}_{DR},\end{equation} where the $\omega$-primitive subspace $H^2_{DR}(X,\,\C)_{\omega\mhyphen prim}$ is a complex hyperplane of $H^2_{DR}(X,\,\C)$ depending only on the cohomology class $\{\omega_{n-1}\}_{DR}\in H^{2n-2}_{DR}(X,\,\C)$, while $\omega_h$ is the $\Delta_\omega$-harmonic component of $\omega$ and the complex line $\C\cdot\{\omega_h\}_{DR}$ depends on the choice of the balanced metric $\omega$.

\end{Conc}

\vspace{2ex}

If $\omega$ is K\"ahler, the Lefschetz-type decomposition (\ref{eqn:Lefschetz-type-decomp_bal}) depends only on the {\it K\"ahler class} $\{\omega\}_{DR}\in H^2_{DR}(X,\,\C)$ since $\omega_h=\omega$ in that case.

\begin{Lem}\label{Lem:Lefschetz-type-decomp_bal_class_coeff} The assumptions are the same as in Conclusion \ref{Conc:Lefschetz-type-decomp_bal}. For every $\alpha\in C^\infty_2(X,\,\C)\cap\ker d$, the coefficient of $\{\omega_h\}_{DR}$ in the Lefschetz-type decomposition of $\{\alpha\}_{DR}\in H^2_{DR}(X,\,\C)$ according to (\ref{eqn:Lefschetz-type-decomp_bal}), namely in \begin{equation}\label{eqn:Lefschetz-type-decomp_bal_class}\{\alpha\}_{DR} = \{\alpha\}_{DR,\,prim} + \lambda\,\{\omega_h\}_{DR},\end{equation} is given by \begin{equation}\label{eqn:Lefschetz-type-decomp_bal_class_coeff}\lambda = \lambda_\omega(\{\alpha\}_{DR}) = \frac{\{\omega_{n-1}\}_{DR}.\{\alpha\}_{DR}}{||\omega_h||^2_\omega} = \frac{1}{||\omega_h||^2_\omega}\,\int\limits_X\alpha\wedge\omega_{n-1}.\end{equation}

\end{Lem}

\noindent {\it Proof.} Since $\{\alpha\}_{DR,\,prim}\in H^2_{DR}(X,\,\C)_{\omega\mhyphen prim}$, we have $\{\omega_{n-1}\}_{DR}.\{\alpha\}_{DR,\,prim}=0$, so \begin{eqnarray*}\{\omega_{n-1}\}_{DR}.\{\alpha\}_{DR} & = & \lambda\,\int\limits\omega_{n-1}\wedge\omega_h = \lambda\,\int\limits(\omega_{n-1})_h\wedge\omega_h = \lambda\,\int\limits(\omega_{n-1})_h\wedge\star_\omega(\omega_{n-1})_h \\
  & = &  \lambda\,||(\omega_{n-1})_h||^2_\omega = \lambda\,||\omega_h||^2_\omega.\end{eqnarray*} This gives (\ref{eqn:Lefschetz-type-decomp_bal_class_coeff}).  \hfill $\Box$

\vspace{2ex}

Formula (\ref{eqn:Lefschetz-type-decomp_bal_class_coeff}) implies that $\lambda_\omega(\{\alpha\}_{DR})$ is {\it real} if the class $\{\alpha\}_{DR}\in H^2_{DR}(X,\,\R)$ is real. This enables one to define a {\it positive side} and a {\it negative side} of the hyperplane $H^2_{DR}(X,\,\R)_{\omega\mhyphen prim}:=H^2_{DR}(X,\,\C)_{\omega\mhyphen prim}\cap H^2_{DR}(X,\,\R)$ in $H^2_{DR}(X,\,\R)$ by \begin{eqnarray*}\label{eqn:positive-negative-sides_DR}\nonumber H^2_{DR}(X,\,\R)_\omega^{+} & := & \bigg\{\{\alpha\}_{DR}\in H^2_{DR}(X,\,\R)\,\mid\,\lambda_\omega(\{\alpha\}_{DR})>0\bigg\}, \\
H^2_{DR}(X,\,\R)_\omega^{-} & := & \bigg\{\{\alpha\}_{DR}\in H^2_{DR}(X,\,\R)\,\mid\,\lambda_\omega(\{\alpha\}_{DR})<0\bigg\}.\end{eqnarray*} These open subsets of $H^2_{DR}(X,\,\R)$ depend only on the cohomology class $\{\omega_{n-1}\}_{DR}\in H^{2n-2}_{DR}(X,\,\R)$.

Since $\{\alpha\}_{DR}$ is $\omega$-primitive if and only if $\lambda_\omega(\{\alpha\}_{DR})=0$, we get a {\it partition} of $H^2_{DR}(X,\,\R)$: \begin{eqnarray*}\label{eq:partition_H2_DR} H^2_{DR}(X,\,\R) = H^2_{DR}(X,\,\R)_\omega^{+}\cup H^2_{DR}(X,\,\R)_{\omega\mhyphen prim}\cup H^2_{DR}(X,\,\R)_\omega^{-}\end{eqnarray*} depending only on the cohomology class $\{\omega_{n-1}\}_{DR}\in H^{2n-2}_{DR}(X,\,\R)$.

\vspace{2ex}

The next (trivial) observation is that the $\omega$-primitive hyperplane $H^2_{DR}(X,\,\C)_{\omega\mhyphen prim}\subset H^2_{DR}(X,\,\C)$ depends only on the ray $\R_{>0}\cdot\{\omega_{n-1}\}_{DR}$ generated by the De Rham cohomology class of $\omega_{n-1}$ in the De Rham version of the {\it balanced cone} ${\cal B}_{X,\,DR}\subset H^{2n-2}_{DR}(X,\,\R)$ of $X$. (We denote by ${\cal B}_{X,\,DR}$ the set of De Rham cohomology classes $\{\omega_{n-1}\}_{DR}$ induced by balanced metrics $\omega$.)

\begin{Lem}\label{Lem:ray_DR_cones} Let $X$ be a compact complex non-degenerate balanced manifold with $\mbox{dim}_\C X=n$. Let $\omega$ and $\gamma$ be {\bf balanced metrics} on $X$. The following equivalence holds: \begin{equation*}H^2_{DR}(X,\,\C)_{\omega\mhyphen prim} = H^2_{DR}(X,\,\C)_{\gamma\mhyphen prim} \iff \exists\,c>0 \hspace{1ex}\mbox{such that}\hspace{1ex} \{\omega_{n-1}\}_{DR} = c\, \{\gamma_{n-1}\}_{DR}.\end{equation*}

\end{Lem}

\noindent {\it Proof.} ``$\Longleftarrow$'' This implication follows from proportional linear maps having the same kernel.

``$\Longrightarrow$'' This implication follows from the following elementary fact. Suppose $T,S:E\longrightarrow\C$ are $\C$-linear maps on a $\C$-vector space $E$ such that $\ker T = \ker S\subset E$ is of $\C$-codimension $1$ in $E$. Then, there exists $c\in\C\setminus\{0\}$ such that $T = cS$. To see this, let $\{e_j\,\mid\, j\in J\}$ be a $\C$-basis of $\ker T = \ker S$ and let $e\in E$ such that $\{e\}\cup\{e_j\,\mid\, j\in J\}$ is a $\C$-basis of $E$. Then, $T(e)$ and $S(e)$ are non-zero complex numbers, so there exists a unique $c\in\C\setminus\{0\}$ such that $T(e) = c\,S(e)$. Now, fix an arbitrary $u\in E$. We will show that $T(u) = c\,S(u)$. There is a unique choice of $\lambda\in\C$ and $v\in\ker T = \ker S$ such that $u=\lambda\, e + v$. Hence, $T(u) = \lambda\,T(e) = c\,(\lambda\,S(e)) = c\,S(u)$.

In our case, the assumption $H^2_{DR}(X,\,\C)_{\omega\mhyphen prim} = H^2_{DR}(X,\,\C)_{\gamma\mhyphen prim}$ amounts to $\ker(\{\omega_{n-1}\}_{DR}\wedge\cdot) = \ker(\{\gamma_{n-1}\}_{DR}\wedge\cdot)$. Hence, by the above elementary fact, it amounts to the exietence of a constant $c\in\C\setminus\{0\}$ such that $\{\omega_{n-1}\}_{DR}\wedge\cdot = c\,\{\gamma_{n-1}\}_{DR}\wedge\cdot$ as $\C$-linear maps on $H^2_{DR}(X,\,\C)$. By the non-degeneracy of the Poincar\'e duality $H^2_{DR}(X,\,\C)\times H^{2n-2}_{DR}(X,\,\C)\longrightarrow\C$, this further amounts to the existence of a constant $c\in\C\setminus\{0\}$ such that $\{\omega_{n-1}\}_{DR} = c\,\{\gamma_{n-1}\}_{DR}$.

Now, since the forms $\omega_{n-1}$ and $\gamma_{n-1}$ are real, the constant $c$ can be chosen real. (Replace $c$ with $(c+\bar{c})/2$ if necessary.) Since the balanced metric $\omega$ is non-degenerate, $c\neq 0$. If $c<0$, then $\omega_{n-1} -c\,\gamma_{n-1}>0$ would be the $d$-exact $(n-1)$-st power of a balanced metric. This balanced metric would then be degenerate balanced, contradicting the assumption on $X$. Thus, $c$ must be positive. \hfill $\Box$

\vspace{2ex}

We will now see that not only do proportional balanced classes $\{\omega_{n-1}\}_{DR}$ and $\{\gamma_{n-1}\}_{DR}$ induce the same hyperplane of primitive classes in $H^2_{DR}(X,\,\C)$, but they can be made to also induce the same Lefschetz-type decomposition (\ref{eqn:Lefschetz-type-decomp_bal}). This is fortunate since, in general, the complex line $\C\cdot\{\omega_h\}_{DR}$ depends on the choice of the balanced metric $\omega$, unlike $H^2_{DR}(X,\,\C)_{\omega\mhyphen prim}$ which depends only on the balanced class $\{\omega_{n-1}\}_{DR}\in H^{2n-2}_{DR}(X,\,\C)$.

\begin{Lem}\label{Lem:ray_DR_cones_Lefschetz-decomp} Let $X$ be a compact complex non-degenerate balanced manifold with $\mbox{dim}_\C X=n$.

  \vspace{1ex}

  (i)\, If $\omega$ and $\gamma$ are {\bf balanced metrics} on $X$ such that $\omega_{n-1} = c\,\gamma_{n-1}$ for some constant $c>0$, there exists a constant $a>0$ such that $\omega_h = a\,\gamma_h$.

  \vspace{1ex}

  (ii)\, For every ray $\R_{>0}\cdot\{\omega_{n-1}\}_{DR}$ in the De Rham version of the {\it balanced cone} ${\cal B}_{X,\,DR}\subset H^{2n-2}_{DR}(X,\,\R)$ of $X$, the balanced metrics representing the classes on this ray can be chosen such that they induce the same Lefschetz-type decomposition (\ref{eqn:Lefschetz-type-decomp_bal}).

\end{Lem}

\noindent{\it Proof.} (i)\, Since $\omega = c^{\frac{1}{n-1}}\,\gamma$, we get $\star_\omega = {const}\cdot\star_\gamma$ and $d^\star_\omega = {const}\cdot d^\star_\gamma$. The latter identity implies $\Delta_\omega = {const}\cdot\Delta_\gamma$, hence $\ker\Delta_\omega = \ker\Delta_\gamma$. In particular, $(\omega_{n-1})_h = c\,(\gamma_{n-1})_h$ and thus $$\omega_h = \star_\omega(\omega_{n-1})_h = {const}\cdot\star_\gamma(\gamma_{n-1})_h = {const}\cdot\gamma_h,$$ where in all the above identities ${const}$ stands for a positive constant that may change from one occurrence to another.

\vspace{1ex}

(ii)\, Fix a balanced De Rham class $\{\gamma_{n-1}\}_{DR}\in{\cal B}_{X,\,DR}\subset H^{2n-2}_{DR}(X,\,\R)$ and fix a balanced metric $\gamma$ (whose choice is arbitrary) such that $\gamma_{n-1}$ represents the class $\{\gamma_{n-1}\}_{DR}$. For every constant $c>0$, the balanced class $c\,\{\gamma_{n-1}\}_{DR}$ can be represented by the form $\omega_{n-1}:=c\,\gamma_{n-1}$ which is induced by the balanced metric $\omega:=c^{\frac{1}{n-1}}\,\gamma$. From (i), we get $\C\,\{\omega_h\}_{DR} = \C\,\{\gamma_h\}_{DR}$. Since we also have $H^2_{DR}(X,\,\C)_{\omega\mhyphen prim} = H^2_{DR}(X,\,\C)_{\gamma\mhyphen prim}$ by Lemma \ref{Lem:ray_DR_cones}, the contention follows.  \hfill $\Box$

\vspace{2ex}

The proof of (ii) of the above Lemma \ref{Lem:ray_DR_cones_Lefschetz-decomp} shows that the line $\C\,\{\omega_h\}_{DR}$ in the Lefschetz-type decomposition (\ref{eqn:Lefschetz-type-decomp_bal}) induced by a given ray $\R_{>0}\cdot\{\omega_{n-1}\}_{DR}$ in the De Rham version of the {\it balanced cone} ${\cal B}_{X,\,DR}\subset H^{2n-2}_{DR}(X,\,\R)$ of $X$ still depends on the arbitrary choice of a balanced metric $\gamma$ such that $\gamma_{n-1}$ represents a given class $\{\gamma_{n-1}\}_{DR}$ on this ray. To tame this dependence, we can fix an arbitrary Hermitian (not necessarily balanced) metric $\rho$ on $X$ and make all the choices of harmonic representatives and projections be induced by $\rho$. Thus, we get $L^2_\rho$-orthogonal decompositions: \begin{equation*}\label{eqn:omega_power_3-space-decomp_rho}\omega = \omega_{h,\,\rho} + d_\rho^\star\,\eta_{\omega,\,\rho} \hspace{3ex} \mbox{and} \hspace{3ex} \omega_{n-1} = (\omega_{n-1})_{h,\,\rho} + d\,\Gamma_{\omega,\,\rho},\end{equation*} where $\omega_{h,\,\rho}\in\ker\Delta_\rho$ as a $2$-form, $(\omega_{n-1})_{h,\,\rho}\in\ker\Delta_\rho$ as a $(2n-2)$-form, while $\eta_{\omega,\,\rho}$ and $\Gamma_{\omega,\,\rho}$ are smooth forms of respective degrees $3$ and $2n-3$. Since $\omega$ and $\omega_{n-1}$ are real, so are their $\Delta_\rho$-harmonic components $\omega_{h,\,\rho}$ and $(\omega_{n-1})_{h,\,\rho}$.

In this way, every {\it non-zero balanced class} $\{\omega_{n-1}\}_{DR}$ induces a Lefschetz-type decomposition analogous to (\ref{eqn:Lefschetz-type-decomp_bal}) that depends only on the class $\{\omega_{n-1}\}_{DR}$ and on the background metric $\rho$: \begin{equation*}\label{eqn:Lefschetz-type-decomp_bal_rho}H^2_{DR}(X,\,\C) = H^2_{DR}(X,\,\C)_{\omega\mhyphen prim}\oplus\C\cdot\{\omega_{h,\,\rho}\}_{DR},\end{equation*} where the hyperplane $H^2_{DR}(X,\,\C)_{\omega\mhyphen prim}$ depends only on the class $\{\omega_{n-1}\}_{DR}$.

In other words, we remove the dependence of the line $\C\,\{\omega_h\}_{DR}$ in the Lefschetz-type decomposition (\ref{eqn:Lefschetz-type-decomp_bal}) on a representative of the class $\{\omega_{n-1}\}_{DR}$ and replace it with the dependence on a fixed background metric $\rho$.

\subsection{Case of degree $2$: Bott-Chern and Aeppli cohomologies}\label{subsection:degree2_BC-A} Let us finally point out that the theory developed in $\S.$\ref{subsection:degree2_DR} in the context of the Poincar\'e duality for the De Rham cohomology spaces of degrees $2$ and $2n-2$ can be rerun in the context of the duality (\ref{eqn:BC-A_duality}) between the Bott-Chern and Aeppli cohomology spaces of bidegrees $(1,\,1)$, resp. $(n-1,\,n-1)$.
 
Since all the results and constructions of $\S.$\ref{subsection:degree2_DR}, except for Lemma \ref{Lem:primitive-classes_rep}, have analogues in the new context with very similar proofs, we will leave most of these proofs to the reader.

In fact, the new cohomological setting allows for the theory of $\S.$\ref{subsection:degree2_DR} to be repeated in the more general context of {\it Gauduchon} (not necessarily balanced) {\it metrics} and the Aeppli cohomology classes they define in $H^{n-1,\,n-1}_A(X,\,\R)$. We start with the following analogue of Lemma and Definition \ref{Lem-Def:primitive-classes_deg2}.

\begin{Lem-Def}\label{Lem-Def:primitive-classes_bideg11} Let $\omega$ be a {\bf Gauduchon} metric on a compact complex manifold $X$ with $\mbox{dim}_\C X=n$. The linear map: \begin{equation*}\label{eqn:primitive-classes_bideg11}[\omega_{n-1}]_A\wedge\cdot : H^{1,\,1}_{BC}(X,\,\C)\longrightarrow H^{n,\,n}_A(X,\,\C)\simeq\C, \hspace{3ex} [\alpha]_{BC}\longmapsto[\omega_{n-1}\wedge\alpha]_A,\end{equation*} is {\bf well defined} and depends only on the cohomology class $[\omega_{n-1}]_A\in H^{n-1,\,n-1}_A(X,\,\C)$. We set: \begin{equation*}H^{1,\,1}_{BC}(X,\,\C)_{\omega\mhyphen prim}:=\ker\bigg([\omega_{n-1}]_A\wedge\cdot\bigg)\subset H^{1,\,1}_{BC}(X,\,\C)\end{equation*} and we call its elements {\bf ($\omega$-)primitive} Bott-Chern $(1,\,1)$-classes.

\end{Lem-Def}

\noindent {\it Proof.} The well-definedness follows at once from the identities: \begin{eqnarray*}\partial\bar\partial(\omega_{n-1}\wedge\alpha) & = & \partial\bar\partial\omega_{n-1}\wedge\alpha = 0,  \hspace{5ex} \alpha\in C^\infty_{1,\,1}(X,\,\C)\cap\ker d,\\
  \omega_{n-1}\wedge\partial\bar\partial\varphi & = & \partial(\omega_{n-1}\wedge\bar\partial\varphi) + \bar\partial(\varphi\,\partial\omega_{n-1})\in\mbox{Im}\,\partial + \mbox{Im}\,\bar\partial, \hspace{5ex} \varphi\in C^\infty_{0,\,0}(X,\,\C),\end{eqnarray*} where the latter takes into account the fact that $\partial\bar\partial\omega_{n-1}=0$.

That the map $[\omega_{n-1}]_A\wedge\cdot$ depends only on the Aeppli cohomology class $[\omega_{n-1}]_A$ follows from: \begin{eqnarray*}(\omega_{n-1} + \partial\bar\Gamma + \bar\partial\Gamma)\wedge\alpha - \omega_{n-1}\wedge\alpha = \partial(\bar\Gamma\wedge\alpha) + \bar\partial(\Gamma\wedge\alpha)\in\mbox{Im}\,\partial + \mbox{Im}\,\bar\partial,\hspace{5ex} \alpha\in C^\infty_{1,\,1}(X,\,\C)\cap\ker d.\end{eqnarray*} \hfill $\Box$

The following result is the analogue of Lemma \ref{Lem:primitive_deg-bal_char}.

\begin{Lem}\label{Lem:primitive_Aeppli-exact_char} Let $\omega$ be a {\bf Gauduchon} metric on a compact complex manifold $X$ with $\mbox{dim}_\C X=n$. The following equivalence holds: $$H^{1,\,1}_{BC}(X,\,\C)_{\omega\mhyphen prim}= H^{1,\,1}_{BC}(X,\,\C) \iff \omega_{n-1}\in\mbox{Im}\,\partial + \mbox{Im}\,\bar\partial \hspace{3ex}  (i.e. \hspace{1ex} \omega_{n-1} \hspace{1ex} \mbox{\bf is Aeppli-exact}).$$

\end{Lem}

\noindent {\it Proof.} ``$\Longleftarrow$'' If $\omega_{n-1}\in\mbox{Im}\,\partial + \mbox{Im}\,\bar\partial$, $[\omega_{n-1}]_A=0$, so the map $[\omega_{n-1}]_A\wedge\cdot$ vanishes identically.

\vspace{1ex}

``$\implies$'' Suppose that $H^{1,\,1}_{BC}(X,\,\C)_{\omega\mhyphen prim}= H^{1,\,1}_{BC}(X,\,\C)$. This translates to \begin{equation}\label{eqn:primitive_Aeppli-exact_char_proof_1}\omega_{n-1}\wedge\alpha\in\mbox{Im}\,\partial + \mbox{Im}\,\bar\partial, \hspace{6ex}\forall\,\alpha\in C^\infty_{1,\,1}(X,\,\C)\cap\ker d.\end{equation}

Since $\omega_{n-1}$ is $(\partial\bar\partial)$-closed, it has a unique $L^2_\omega$-orthogonal decomposition: \begin{equation*}\omega_{n-1} = (\omega_{n-1})_h + (\partial\bar\Gamma_\omega + \bar\partial\Gamma_\omega),\end{equation*} with an $(n-1,\,n-1)$-form $(\omega_{n-1})_h\in\ker\Delta_{A,\,\omega}$ and an $(n-1,\,n-2)$-form $\Gamma_\omega$. (See (\ref{eqn:BC-A_3-space-decomp}) below.)

On the other hand, $\alpha$ is $d$-closed, so it has a unique $L^2_\omega$-orthogonal decomposition: \begin{equation*}\alpha = \alpha_h + \partial\bar\partial\varphi,\end{equation*} where $\alpha_h$ is $\Delta_{BC,\,\omega}$-harmonic and $\varphi$ is a smooth function on $X$. (See again (\ref{eqn:BC-A_3-space-decomp}) below.)

Thus, for every $\alpha\in C^\infty_{1,\,1}(X,\,\C)\cap\ker d$, we get: \begin{eqnarray*}\omega_{n-1}\wedge\alpha & = & (\omega_{n-1})_h\wedge\alpha + \partial(\bar\Gamma_\omega\wedge\alpha) + \bar\partial(\Gamma_\omega\wedge\alpha) \\
  & = & (\omega_{n-1})_h\wedge\alpha_h + \partial\bigg((\omega_{n-1})_h\wedge\bar\partial\varphi\bigg) + \bar\partial\bigg(\varphi\,\partial(\omega_{n-1})_h\bigg) + \partial(\bar\Gamma_\omega\wedge\alpha) + \bar\partial(\Gamma_\omega\wedge\alpha),\end{eqnarray*} where for the last identity we used the fact that $\bar\partial\partial(\omega_{n-1})_h=0$.

Thanks to assumption (\ref{eqn:primitive_Aeppli-exact_char_proof_1}), the last identity implies that \begin{equation}\label{eqn:primitive_Aeppli-exact_char_proof_2}(\omega_{n-1})_h\wedge\alpha_h\in\mbox{Im}\,\partial + \mbox{Im}\,\bar\partial, \hspace{6ex}\forall\,\alpha_h\in C^\infty_{1,\,1}(X,\,\C)\cap\ker\Delta_{BC,\,\omega}.\end{equation}

Meanwhile, since $(\omega_{n-1})_h$ is $\Delta_{A,\,\omega}$-harmonic (and real), $\star_\omega(\omega_{n-1})_h$ is $\Delta_{BC,\,\omega}$-harmonic (and real). Hence, \begin{equation*}\mbox{Im}\,\partial + \mbox{Im}\,\bar\partial\ni(\omega_{n-1})_h\wedge\star_\omega(\omega_{n-1})_h = |(\omega_{n-1})_h|^2_\omega\,dV_\omega\geq 0,\end{equation*} where the first relation follows from (\ref{eqn:primitive_Aeppli-exact_char_proof_2}) by choosing $\alpha_h = \star_\omega(\omega_{n-1})_h$. Consequently, from Stokes's Theorem we get: $$\int\limits_X |(\omega_{n-1})_h|^2_\omega\,dV_\omega = 0,$$ hence $(\omega_{n-1})_h=0$. This implies that $\omega_{n-1}\in\mbox{Im}\,\partial + \mbox{Im}\,\bar\partial$ and we are done.   \hfill $\Box$

\vspace{2ex}

The analogue in this context of Corollary \ref{Cor:primitive_bal_deg-bal} is the following

\begin{Cor}\label{Cor:primitive_Gauduchon_Aeppli-exact} Let $\omega$ be a {\bf Gauduchon} metric on a compact complex manifold $X$ with $\mbox{dim}_\C X=n$. The following dichotomy holds:

  \vspace{1ex}

  (a)\, if $\omega_{n-1}$ is {\bf not Aeppli exact}, $H^{1,\,1}_{BC}(X,\,\C)_{\omega\mhyphen prim}$ is a complex hyperplane of $H^{1,\,1}_{BC}(X,\,\C)$ depending only on the Aeppli-Gauduchon class $\{\omega_{n-1}\}_{DR}\in{\cal G}$;

  \vspace{1ex}

  (b)\, if $\omega_{n-1}$ is {\bf Aeppli exact}, $H^{1,\,1}_{BC}(X,\,\C)_{\omega\mhyphen prim}= H^{1,\,1}_{BC}(X,\,\C)$.
  
\end{Cor}

\vspace{2ex}

To get a {\it Lefschetz-type decomposition} of $H^{1,\,1}_{BC}(X,\,\C)$ induced by an arbitrary Gauduchon metric $\omega$, we use the orthogonal $3$-space decompositions featuring the Aeppli-, resp. Bott-Chern-Laplacians induced by the metric $\omega$: \begin{eqnarray}\label{eqn:BC-A_3-space-decomp}\nonumber C^\infty_{n-1,\,n-1}(X,\,\C) & = & \ker\Delta_{A,\,\omega}\oplus(\mbox{Im}\,\partial + \mbox{Im}\,\bar\partial)\oplus\mbox{Im}\,(\partial\bar\partial)^\star, \\
  C^\infty_{1,\,1}(X,\,\C) & = & \ker\Delta_{BC,\,\omega}\oplus\mbox{Im}\,(\partial\bar\partial)\oplus(\mbox{Im}\,\partial^\star + \mbox{Im}\,\bar\partial^\star),\end{eqnarray} where $\ker\Delta_{A,\,\omega}\oplus(\mbox{Im}\,\partial + \mbox{Im}\,\bar\partial) = \ker (\partial\bar\partial)$ and $\ker\Delta_{BC,\,\omega}\oplus(\mbox{Im}\,\partial^\star + \mbox{Im}\,\bar\partial^\star) = \ker(\partial\bar\partial)^\star$. Thus, we get unique decompositions of $\omega$, resp. $\omega_{n-1}$: \begin{equation}\label{eqn:omega_power_3-space-decomp_BC-A}\ker(\partial\bar\partial)^\star\ni\omega = \omega_h + (\partial_\omega^\star\bar{u}_\omega + \bar\partial_\omega^\star u_\omega) \hspace{3ex} \mbox{and} \hspace{3ex} \ker(\partial\bar\partial)\ni\omega_{n-1} = (\omega_{n-1})_h + (\partial\bar\Gamma_\omega + \bar\partial\Gamma_\omega),\end{equation} where $\omega_h\in\ker\Delta_{BC,\,\omega}$ as a $(1,\,1)$-form, $(\omega_{n-1})_h\in\ker\Delta_{A,\,\omega}$ as an $(n-1,\,n-1)$-form, while $u_\omega$ and $\Gamma_\omega$ are smooth forms of respective bidegrees $(1,\,2)$ and $(n-1,\,n-2)$. Since $\omega$ and $\omega_{n-1}$ are real, so are their harmonic components $\omega_h$ and $(\omega_{n-1})_h$. Since $\star_\omega\omega = \omega_{n-1}$ and since the Hodge star operator $\star_\omega$ maps Aeppli-harmonic forms to Bott-Chern-harmonic forms and vice-versa, we get: \begin{equation}\label{eqn:omega_power_3-star-interchange_BC-A}\star_\omega\omega_h = (\omega_{n-1})_h \hspace{3ex} \mbox{and} \hspace{3ex} \star_\omega(\partial_\omega^\star\bar{u}_\omega + \bar\partial_\omega^\star u_\omega) = \partial\bar\Gamma_\omega + \bar\partial\Gamma_\omega.\end{equation} 

Thus, $\omega_h$ is uniquely determined by $\omega$ and is $d$-closed (because it is even $\Delta_{BC,\,\omega}$-harmonic). Therefore, it represents a class in $H^{1,\,1}_{BC}(X,\,\R)$.

\begin{Def}\label{Def:H11-BC-class_Gauduchon-metric} For any Gauduchon metric $\omega$ on a compact complex manifold $X$, the Bott-Chern cohomlogy class $[\omega_h]_{BC}\in H^{1,\,1}_{BC}(X,\,\R)$ is called the {\bf cohomology class} of $\omega$.

\end{Def}  

Of course, if $\omega$ is K\"ahler, $\omega_h = \omega$, so $\{\omega_h\}_{BC}$ is the usual Bott-Chern K\"ahler class $\{\omega\}_{BC}$.

\vspace{1ex}

The analogue of Lemma \ref{Lem:orthogonality_prim-class-rep_balanced} is the following

\begin{Lem}\label{Lem:orthogonality_prim-class-rep_Gauduchon} Suppose there exists a Gauduchon metric $\omega$ on a compact complex manifold $X$. Then, for every $\alpha\in C^\infty_{1,\,1}(X,\,\C)$ such that $d\alpha=0$ and $[\alpha]_{BC}\in H^{1,\,1}_{BC}(X,\,\C)_{\omega\mhyphen prim}$, we have: $$\langle\langle\omega_h,\,\alpha\rangle\rangle_\omega = 0,$$ where $\langle\langle\,\,,\,\,\rangle\rangle_\omega$ is the $L^2$ inner product induced by $\omega$.

\end{Lem}

%\noindent {\it Proof.} It is virtually identical to the proof of Lemma \ref{Lem:orthogonality_prim-class-rep_balanced}. \hfill $\Box$

\vspace{2ex} The analogue in this context of Conclusion \ref{Conc:Lefschetz-type-decomp_bal} is the following

\begin{Conc}\label{Conc:Lefschetz-type-decomp_Gauduchon} Let $X$ be a compact complex manifold with $\mbox{dim}_\C X=n$. Let $\omega$ be a {\bf Gauduchon metric} on $X$ such that $\omega_{n-1}$ is not Aeppli-exact. Then, the Bott-Chern cohomology space of bidegree $(1,\,1)$ has a {\bf Lefschetz-type} $L^2_\omega$-orthogonal {\bf decomposition}: \begin{equation}\label{eqn:Lefschetz-type-decomp_Gauduchon}H^{1,\,1}_{BC}(X,\,\C) = H^{1,\,1}_{BC}(X,\,\C)_{\omega\mhyphen prim}\oplus\C\cdot[\omega_h]_{BC},\end{equation} where the $\omega$-primitive subspace $H^{1,\,1}_{BC}(X,\,\C)_{\omega\mhyphen prim}$ is a complex hyperplane of $H^{1,\,1}_{BC}(X,\,\C)$ depending only on the cohomology class $[\omega_{n-1}]_A\in H^{n-1,\,n-1}_A(X,\,\C)$, while $\omega_h$ is the $\Delta_{BC,\,\omega}$-harmonic component of $\omega$ and the complex line $\C\cdot[\omega_h]_{BC}$ depends on the choice of the Gauduchon metric $\omega$.

\end{Conc}

We also have the following analogue of Lemma \ref{Lem:Lefschetz-type-decomp_bal_class_coeff}.

\begin{Lem}\label{Lem:Lefschetz-type-decomp_Gauduchon_class_coeff} The assumptions are the same as in Conclusion \ref{Conc:Lefschetz-type-decomp_Gauduchon}. For every $\alpha\in C^\infty_{1,\,1}(X,\,\C)\cap\ker d$, the coefficient of $[\omega_h]_{BC}$ in the Lefschetz-type decomposition of $[\alpha]_{BC}\in H^{1,\,1}_{BC}(X,\,\C)$ according to (\ref{eqn:Lefschetz-type-decomp_Gauduchon}), namely in \begin{equation}\label{eqn:Lefschetz-type-decomp_Aeppli-Gauduchon_class}[\alpha]_{BC} = [\alpha]_{BC,\,prim} + \lambda\,[\omega_h]_{BC},\end{equation} is given by \begin{equation}\label{eqn:Lefschetz-type-decomp_Aeppli-Gauduchon_class_coeff}\lambda = \lambda_\omega([\alpha]_{BC}) = \frac{[\omega_{n-1}]_A.[\alpha]_{BC}}{||\omega_h||^2_\omega} = \frac{1}{||\omega_h||^2_\omega}\,\int\limits_X\alpha\wedge\omega_{n-1}.\end{equation}

\end{Lem}

%\noindent {\it Proof.} It is virtually identical to the proof of Lemma \ref{Lem:Lefschetz-type-decomp_bal_class_coeff}.   \hfill $\Box$

As in $\S.$\ref{subsection:degree2_DR}, formula (\ref{eqn:Lefschetz-type-decomp_Aeppli-Gauduchon_class_coeff}) implies that $\lambda_\omega([\alpha]_{BC})$ is {\it real} if the class $[\alpha]_{BC}\in H^{1,\,1}_{BC}(X,\,\C)$ is real. Thus, we can define a {\it positive side} and a {\it negative side} of the hyperplane $H^{1,\,1}_{BC}(X,\,\R)_{\omega\mhyphen prim}:=H^{1,\,1}_{BC}(X,\,\C)_{\omega\mhyphen prim}\cap H^{1,\,1}_{BC}(X,\,\R)$ in $H^{1,\,1}_{BC}(X,\,\R)$ by \begin{eqnarray*}\label{eqn:positive-negative-sides_BC}\nonumber H^{1,\,1}_{BC}(X,\,\R)_\omega^{+} & := & \bigg\{[\alpha]_{BC}\in H^{1,\,1}_{BC}(X,\,\R)\,\mid\,\lambda_\omega([\alpha]_{BC})>0\bigg\}, \\
H^{1,\,1}_{BC}(X,\,\R)_\omega^{-} & := & \bigg\{[\alpha]_{BC}\in H^{1,\,1}_{BC}(X,\,\R)\,\mid\,\lambda_\omega([\alpha]_{BC})<0\bigg\}.\end{eqnarray*} These are open subsets of $H^{1,\,1}_{BC}(X,\,\R)$ that depend only on the cohomology class $[\omega_{n-1}]_A\in H^{n-1,\,n-1}_A(X,\,\R)$.

Since $[\alpha]_{BC}$ is $\omega$-primitive if and only if $\lambda_\omega([\alpha]_{BC})=0$, we get a {\it partition} of $H^{1,\,1}_{BC}(X,\,\R)$: \begin{eqnarray*}\label{eq:partition_H2_BC} H^{1,\,1}_{BC}(X,\,\R) = H^{1,\,1}_{BC}(X,\,\R)_\omega^{+}\cup H^{1,\,1}_{BC}(X,\,\R)_{\omega\mhyphen prim}\cup H^{1,\,1}_{BC}(X,\,\R)_\omega^{-}\end{eqnarray*} depending only on the cohomology class $[\omega_{n-1}]_A\in H^{n-1,\,n-1}_A(X,\,\R)$.

\vspace{2ex}

As a consequence of these considerations, we get

\begin{Prop}\label{Prop:psef-cone_positive-side-char} Let $X$ be a compact complex manifold with $\mbox{dim}_\C X=n$. The {\bf pseudo-effective cone} ${\cal E}_X\subset H^{1,\,1}_{BC}(X,\,\R)$ of $X$ is the intersection of all the {\bf non-negative sides} $$H^{1,\,1}_{BC}(X,\,\R)_\omega^{\geq 0}:=H^{1,\,1}_{BC}(X,\,\R)_\omega^{+}\cup H^{1,\,1}_{BC}(X,\,\R)_{\omega\mhyphen prim}$$  of hyperplanes $H^{1,\,1}_{BC}(X,\,\R)_{\omega\mhyphen prim}$ determined by Aeppli-Gauduchon classes $[\omega_{n-1}]_A\in{\cal G}_X$: \begin{equation}\label{eqn:psef-cone_positive-side-char}{\cal E}_X = \bigcap\limits_{[\omega_{n-1}]_A\in{\cal G}_X} H^{1,\,1}_{BC}(X,\,\R)_\omega^{\geq 0},\end{equation}

\end{Prop}

\noindent {\it Proof.} By the duality between the pseudo-effective cone ${\cal E}_X$ and the closure $\overline{\cal G}_X$ of the Gauduchon cone (see $\S.$\ref{subsection:deg-bal_background}), we know that a given class $[T]_{BC}\in H^{1,\,1}_{BC}(X,\,\R)$ lies in ${\cal E}_X$ (i.e. $[T]_{BC}$ can be represented by a closed {\it semi-positive} $(1,\,1)$-current) if and only if $$\int\limits_X T\wedge\omega_{n-1}\geq 0 \hspace{5ex} \mbox{for all}\hspace{1ex} [\omega_{n-1}]_A\in{\cal G}_X.$$ The last condition is equivalent to $\lambda_\omega([T]_{BC})\geq 0$, hence to $[T]_{BC}\in H^{1,\,1}_{BC}(X,\,\R)_\omega^{\geq 0}$, for all $[\omega_{n-1}]_A\in{\cal G}_X$, so the contention follows. \hfill $\Box$

\vspace{2ex}

Based on these considerations, we propose Question \ref{Question:small_psef-cone} as a problem for further study.

%\begin{Question} Is it true that the pseudo-effective cone ${\cal E}_X$ is {\bf small} (in a sense to be determined) if (and only if) $X$ is {\bf balanced hyperbolic}?

%\end{Question}  

\section{Properties of balanced hyperbolic manifolds}\label{section:bal-hyperbolic_prop}

The discussion of balanced hyperbolic manifolds featured in this section will mirror that of degenerate balanced manifolds of the previous section.

\subsection{Background and $L^1$ currents on the universal cover}\label{subsection:L1-currents_universal-cover}

It is a classical fact due to Gaffney [Gaf54] that certain basic facts in the Hodge Theory of compact Riemannian manifolds remain valid on {\it complete} such manifolds. The main ingredient in the proof of this fact is the following {\it cut-off trick} of Gaffney's that played a key role in [Gro91, $\S.1$]. It also appears in [Dem97, VIII, Lemma 2.4].

\begin{Lem}([Gaf54])\label{Lem:complete-manifolds_cut-off_functions} Let $(X,\,g)$ be a Riemannian manifold. Then, $(X,\,g)$ is {\bf complete} if and only if there exists an exhaustive sequence $(K_\nu)_{\nu\in\N}$ of {\bf compact} subsets of $X$: $$K_\nu\subset\mathring{K}_{\nu+1} \hspace{3ex} \mbox{for all}\hspace{1ex}\nu\in\N \hspace{3ex} \mbox{and}\hspace{3ex} X = \bigcup\limits_{\nu\in\N}K_\nu,$$ and a sequence $(\psi_\nu)_{\nu\in\N}$ of $C^\infty$ functions $\psi_\nu:X\longrightarrow[0,\,1]$ satisfying, for every $\nu\in\N$, the conditions: \\

\hspace{6ex}  $\psi_\nu = 1$ in a neighbourhood of $K_\nu$,  \hspace{3ex} $\mbox{Supp}\,\psi_\nu\subset\mathring{K}_{\nu+1}$ \hspace{3ex} and

  \vspace{2ex}

  \hspace{6ex}  $||d\psi_\nu||_{L^\infty_g}:=\sup\limits_{x\in X}|(d\psi_\nu)(x)|_g\leq\varepsilon_\nu$,\\

\noindent for some constants $\varepsilon_\nu>0$ such that $\varepsilon_\nu\downarrow 0$ as $\nu$ tends to $+\infty$.    

\end{Lem}  

In particular, the {\it cut-off functions} $\psi_\nu$ are {\it compactly supported}. One can choose $\varepsilon_\nu = 2^{-\nu}$ for each $\nu$ (see e.g. [Dem97, VIII, Lemma 2.4]), but this will play no role here.

\vspace{2ex}

An immediate consequence of Gaffney's cut-off trick is the following classical generalisation of Stokes's Theorem to possibly non-compact, but complete Riemannian manifolds when the forms involved are $L^1$. 

\begin{Lem}([Gro91, Lemma 1.1.A.])\label{Lem:Stokes_non-compact_L1} Let $(X,\,g)$ be a {\bf complete} Riemannian manifold of real dimension $m$. Let $\eta$ be an $L^1_g$-form on $X$ of degree $m-1$ such that $d\eta$ is again $L^1_g$. Then $$\int_Xd\eta = 0.$$

\end{Lem}  

By the form $\eta$ being $L^1$ with respect to the Riemannian metric $g$ ($L^1_g$ for short) we mean that its $L^1$-norm is {\it finite}: $$||\eta||_{L^1_g}:= \int_X|\eta(x)|_g\,dV_g(x)< +\infty,$$ where $dV_g$ is the volume form induced by $g$.

\vspace{2ex}

\noindent {\it Proof of Lemma \ref{Lem:Stokes_non-compact_L1}.} Let $(\psi_\nu)_{\nu\in\N}$ be a sequence of cut-off functions as in Lemma \ref{Lem:complete-manifolds_cut-off_functions} whose existence is guaranteed by the completeness of $(X,\,g)$. The $(m-1)$-form $\psi_\nu\eta$ is {\it compactly supported} for every $\nu\in\N^\star$, so the usual Stokes's Theorem yields: $$\int_Xd(\psi_\nu\eta) = 0,  \hspace{6ex} \nu\in\N^\star.$$

Meanwhile, $d(\psi_\nu\eta) = d\psi_\nu\wedge\eta + \psi_\nu\, d\eta$, so we get: \begin{equation}\label{eqn:Stokes_non-compact_L1_proof_1}\bigg|\int_X\psi_\nu\, d\eta\bigg| = \bigg|\int_X d\psi_\nu\wedge\eta\bigg|\leq||d\psi_\nu||_{L^\infty_g}\,||\eta||_{L^1_g}\leq\varepsilon_\nu\,||\eta||_{L^1_g}, \hspace{6ex} \nu\in\N,\end{equation} for some sequence of constants $\varepsilon_\nu\downarrow 0$.

Since $\eta$ is $L^1_g$, $\varepsilon_\nu\,||\eta||_{L^1_g}\downarrow 0$ as $\nu\to +\infty$. On the other hand, since $d\eta$ is $L^1_g$, the properties of the functions $\psi_\nu$ imply that $$\lim\limits_{\nu\to +\infty}\int_X \psi_\nu d\eta =  \int_X d\eta.$$ Together with (\ref{eqn:Stokes_non-compact_L1_proof_1}), these arguments yield $\int_Xd\eta = 0$, as desired.  \hfill $\Box$

\vspace{2ex}

We now apply this standard cut-off function technique to prove Proposition \ref{Prop:no-L1-currents} stated in the introduction. It is an analogue in our more general context of {\it balanced hyperbolic} manifolds of Proposition 5.4 in [Pop15] according to which compact {\it degenerate balanced} manifolds are characterised by the absence of non-zero $d$-closed positive $(1,\,1)$-currents.

Note that, due to $X$ being {\it compact}, any pair of Hermitian metrics $\omega_1$ and $\omega_2$ on $X$ are {\it comparable} in the sense that there exists a constant $C>0$ such that $(1/C)\,\omega_1\leq\omega_2\leq C\,\omega_1$. Thus, their lifts $\widetilde\omega_1:=\pi^\star\omega_1$ and $\widetilde\omega_2:=\pi^\star\omega_2$ are again comparable on $\widetilde{X}$ by means of the same constant: $(1/C)\,\widetilde\omega_1\leq\widetilde\omega_2\leq C\,\widetilde\omega_1$. Therefore, the $L^1_{\widetilde\omega}$-assumption on $\widetilde{T}$ is independent of the choice of Hermitian metric on $\widetilde{X}$ if this metric is obtained by lifting a metric on $X$. However, the $L^1$-condition changes for metrics on $\widetilde{X}$ that are not lifts of metrics on $X$. But we will not deal with the latter type of metrics.

\vspace{2ex}

\noindent {\it Proof of Proposition \ref{Prop:no-L1-currents}.} Let $n=\mbox{dim}_\C X$. The balanced hyperbolic assumption on $X$ means that $\pi^\star\omega_{n-1} = d\widetilde\Gamma$ on $\widetilde{X}$ for some smooth $L^\infty_{\widetilde\omega}$-form $\widetilde\Gamma$ of degree $(2n-3)$ on $\widetilde{X}$.

If a current $\widetilde{T}$ as in the statement existed on $\widetilde{X}$, we would have \begin{equation}\label{eqn:The:no-L1-currents_proof_1}0<\int\limits_{\widetilde{X}}\widetilde{T}\wedge\pi^\star\omega_{n-1} = \int\limits_{\widetilde{X}}d(\widetilde{T}\wedge\widetilde\Gamma) = 0,\end{equation} which is contradictory.

The last identity in (\ref{eqn:The:no-L1-currents_proof_1}) follows from Lemma \ref{Lem:Stokes_non-compact_L1} applied on the complete manifold $(\widetilde{X},\, \widetilde\omega)$ to the $L^1_{\widetilde\omega}$-current $\eta:= \widetilde{T}\wedge\widetilde\Gamma$ of degree $2n-1$ whose differential $d\eta = \widetilde{T}\wedge\pi^\star\omega_{n-1}$ is again $L^1_{\widetilde\omega}$. That $\eta$ is $L^1_{\widetilde\omega}$ follows from $\widetilde{T}$ being $L^1_{\widetilde\omega}$ (by hypothesis) and from $\widetilde\Gamma$ being $L^\infty_{\widetilde\omega}$, while $d\eta$ being $L^1_{\widetilde\omega}$ follows from $\widetilde{T}$ being $L^1_{\widetilde\omega}$ and from $\pi^\star\omega_{n-1}$ being $L^\infty_{\widetilde\omega}$ (as a lift of the smooth, hence bounded, form $\omega_{n-1}$ on the {\it compact} manifold $X$).  \hfill $\Box$

\vspace{2ex}

We now recall the following standard result saying that some further key facts in the Hodge Theory of compact Riemannian manifolds remain valid on {\it complete} such manifolds $X$ when the differential operators involved (e.g. $d$, $d^\star$, $\Delta$) are considered as {\it closed} and {\it densely defined} unbounded operators on the spaces $L^2_k(X,\,\C)$ of $L^2$-forms of degree $k$ on $X$. The only major property that is lost in passing to complete manifolds is the closedness of the images of these operators. As usual, any differential operator $P$ originally defined on $C^\infty_{\bullet}(X,\,\C)$ is extended to an unbounded operator on $L^2_{\bullet}(X,\,\C)$ by defining its domain $\mbox{Dom}\,P$ as the space of $L^2$-forms $u$ such that $Pu$, computed in the sense of distributions, is again $L^2$.

\begin{The}\label{The:standard_complete_operators}(see e.g. [Dem97, VIII, Theorem 3.2.]) Let $(X,\,g)$ be a {\bf complete} Riemannian manifold of real dimension $m$. Then:

  \vspace{1ex}

(a)\, The space ${\cal D}_{\bullet}(X,\,\C)$ of compactly supported $C^\infty$ forms of any degree (indicated by a $\bullet$) on $X$ is {\bf dense} in the domains $\mbox{Dom}\,d$, $\mbox{Dom}\,d^\star$ and in $\mbox{Dom}\,d\cap\mbox{Dom}\,d^\star$ for the respective graph norms: $$u\mapsto||u|| + ||du||, \hspace{3ex} u\mapsto||u|| + ||d^\star u||, \hspace{3ex} u\mapsto||u|| + ||du|| + ||d^\star u||.$$    

\vspace{1ex}

(b)\, The extension $d^\star$ of the formal adjoint of $d$ to the $L^2$-space coincides with the Hilbert space adjoint of the extension of $d$.  

\vspace{1ex}

(c)\, The $d$-Laplacian $\Delta= \Delta_g:=dd^\star + d^\star d$ has the following property : \begin{equation}\label{eqn:standard_complete_operators_1}\langle\langle\Delta u,\,u\rangle\rangle = ||du||^2 + ||d^\star u||^2\end{equation} for every form $u\in\mbox{Dom}\,\Delta$. In particular, $\mbox{Dom}\,\Delta\subset\mbox{Dom}\,d\cap\mbox{Dom}\,d^\star$ and $\ker\Delta = \ker d\cap\ker d^\star$.

\vspace{1ex}

(d)\, There are {\bf $L^2$-orthogonal decompositions} in every degree (indicated by a $\bullet$): \begin{eqnarray}\label{eqn:standard_complete_operators_2}\nonumber L^2_{\bullet}(X,\,\C) & = & {\cal H}^\bullet_\Delta(X,\,\C)\oplus\overline{\mbox{Im}\,d}\oplus\overline{\mbox{Im}\,d^\star} \\
  \ker d & = & {\cal H}^\bullet_\Delta(X,\,\C)\oplus\overline{\mbox{Im}\,d} \hspace{3ex} \mbox{and} \hspace{3ex} \ker d^\star = {\cal H}^\bullet_\Delta(X,\,\C)\oplus \overline{\mbox{Im}\,d^\star},\end{eqnarray} where ${\cal H}^\bullet_\Delta(X,\,\C):=\{u\in L^2_{\bullet}(X,\,\C)\,\mid\,\Delta u=0\}$ is the space of $\Delta$-{\bf harmonic} $L^2$-forms, while $$\mbox{Im}\,d:= L^2_{\bullet}(X,\,\C)\cap d(L^2_{\bullet -1}(X,\,\C)) \hspace{3ex} \mbox{and} \hspace{3ex} \mbox{Im}\,d^\star:= L^2_{\bullet}(X,\,\C)\cap d^\star(L^2_{\bullet +1}(X,\,\C)).$$

\end{The}

\vspace{2ex}

An immediate consequence of (\ref{eqn:standard_complete_operators_1}) applied in degree $0$ is that on a connected {\it complete} Riemannian manifold $(X,\,g)$, every $\Delta$-{\bf harmonic} $L^2$-function is {\bf constant}: \begin{eqnarray}\label{eqn:harmonic_L2-functions_complete_constant}{\cal H}^0_\Delta(X,\,\C)\subset\C.\end{eqnarray}

\subsection{Harmonic $L^2$-forms of degree $1$ on the universal cover of a balanced hyperbolic manifold}\label{subsection:harmonic_bal-hyp_degree1}

Let $X$ be a possibly non-compact complex manifold with $\mbox{dim}_\C X=n$, supposed to carry a {\it complete balanced metric} $\omega$. In subsequent applications, the roles of $X$ and $\omega$ will be played by $\widetilde{X}$, the universal cover $\pi:\widetilde{X}\longrightarrow X$ of a compact balanced hyperbolic manifold $(X,\,\omega)$, resp. $\widetilde\omega:=\pi^\star\omega$.

A well-known consequence of the K\"ahler commutation relations is the fact that, if $\omega$ is {\it K\"ahler}, the induced $d$-Laplacian $\Delta=\Delta_\omega$ commutes with the multiplication operator $\omega^l\wedge\cdot$ acting on differential forms of any degree on $X$, for every $l$.

We will see that, when $\omega$ is merely {\it balanced}, the commutation of $\Delta$ with the multiplication operator $\omega^{n-1}\wedge\cdot$ acting on differential forms no longer holds. However, we will now compute this commutation defect on $1$-forms.

The computation will continue that of (i) in Lemma \ref{Lem:1-forms_Delta-harm}. For the sake of generality and for a reason that will become apparent later on, we will work with the more general operators $$d_h:=h\partial + \bar\partial,   \hspace{5ex} h\in\C^\star,$$ acting on $\C$-valued differential forms on $X$ and the associated Laplacians $\Delta_h:=d_hd_h^\star + d_h^\star d_h$.

The first stages of the computation lead to the following result in which no completeness assumption is necessary.

\begin{Lem}\label{Lem:commutation-defect_1-forms_pre-integration} Let $X$ be a complex manifold with $\mbox{dim}_\C X=n$. Suppose there exists a {\bf balanced metric} $\omega$ on $X$. Then, for any $h\in\C^\star$ and any {\bf $1$-form} $\varphi$ on $X$, the following identity holds: \begin{eqnarray}\label{eqn:commutation-defect_1-forms_pre-integration}[\Delta_h,\,L_{\omega_{n-1}}]\varphi = \bigg(|h|^2d_{-\frac{1}{\bar{h}}}d_{-\frac{1}{\bar{h}}}^\star - d_h^\star d_h\bigg)\varphi\wedge\omega_{n-1} -i\bar{h}\,d_{-\frac{1}{\bar{h}}}\varphi\wedge d_h\omega_{n-2} -i(|h|^2 +1)\,\partial\bar\partial\varphi\wedge\omega_{n-2}.\end{eqnarray}

\end{Lem}  

\noindent {\it Proof.} $\bullet$ The Jacobi identity yields: $$[[d_h,\,d_h^\star],\,L_{\omega_{n-1}}] - [[d_h^\star,\,L_{\omega_{n-1}}],\,d_h] + [[L_{\omega_{n-1}},\,d_h],\,d_h^\star] = 0.$$ Since $\omega$ is {\it balanced}, $[L_{\omega_{n-1}},\,d_h]=0$. Writing $\Delta_h= [d_h,\,d_h^\star]$, the above equality reduces to \begin{eqnarray}\label{eqn:Delta_h_omega_n-1_Jacobi}[\Delta_h,\,L_{\omega_{n-1}}] = [d_h^\star,\,L_{\omega_{n-1}}]\,d_h + d_h\,[d_h^\star,\,L_{\omega_{n-1}}].\end{eqnarray}

\vspace{1ex}

$\bullet$ Note also the following formula for the {\bf formal adjoint} of $d_h$ involving the Hodge star operator: \begin{eqnarray}\label{eqn:d_h-adjoint_star-formula}d_h^\star = -\bar{h}\star d_{\frac{1}{\bar{h}}}\star.\end{eqnarray}

Indeed, $d_h^\star = (h\partial + \bar\partial)^\star = \bar{h}\,(-\star\bar\partial\star) + (-\star\partial\star) = -\bar{h}\,\star(\frac{1}{\bar{h}}\,\partial + \bar\partial)\star = -\bar{h}\star d_{\frac{1}{\bar{h}}}\star$. No assumption on $\omega$ is needed here.

\vspace{1ex}

$\bullet$ As an application of (\ref{eqn:d_h-adjoint_star-formula}), we observe the following formula for every $(1,\,1)$-form $\alpha$: \begin{eqnarray}\label{eqn:d_h_star_omegan-1_11-form}d_h^\star(\omega_{n-1}\wedge\alpha) = -i\bar{h}\,d_{-\frac{1}{\bar{h}}}(\Lambda_\omega\alpha)\wedge\omega_{n-1}.\end{eqnarray} Again, no assumption on $\omega$ is needed.

To see this, we first multiply the Lefschetz decomposition (\ref{eqn:Lefschetz_decomp_11}) of $\alpha$ by $\omega_{n-1}$ and we get: $\omega_{n-1}\wedge\alpha = (\Lambda_\omega\alpha)\,\omega_n$. Hence, $\star(\omega_{n-1}\wedge\alpha) = \Lambda_\omega\alpha$, so we get the first equality below: $$-\bar{h}\,\star d_{\frac{1}{\bar{h}}}\star(\omega_{n-1}\wedge\alpha) = -\bar{h}\,\star d_{\frac{1}{\bar{h}}}(\Lambda_\omega\alpha) = -\bar{h}\,\star\bigg(\frac{1}{\bar{h}}\,\partial(\Lambda_\omega\alpha)\bigg) - \bar{h}\,\star\bar\partial(\Lambda_\omega\alpha).$$ Applying (\ref{eqn:d_h-adjoint_star-formula}) to the l.h.s. term above and the standard formula (\ref{eqn:prim-form-star-formula-gen}) to the r.h.s. term, we get: $$d_h^\star(\omega_{n-1}\wedge\alpha) = i\partial(\Lambda_\omega\alpha)\wedge\omega_{n-1} -i\bar{h}\,\bar\partial(\Lambda_\omega\alpha)\wedge\omega_{n-1}.$$ Since $i\partial-i\bar{h}\,\bar\partial = -i\bar{h}\,d_{-\frac{1}{\bar{h}}}$, the above equality is nothing but (\ref{eqn:d_h_star_omegan-1_11-form}).

\vspace{1ex}

$\bullet$ {\it Computation of the first term on the r.h.s. of (\ref{eqn:Delta_h_omega_n-1_Jacobi}) on $1$-forms $\varphi = \varphi^{1,\,0} + \varphi^{0,\,1}$.}

\vspace{1ex}

Using formula (\ref{eqn:d_h_star_omegan-1_11-form}) with $\alpha:=h\partial\varphi^{0,\,1} + \bar\partial\varphi^{1,\,0}$, we get the second equality below: \begin{eqnarray}\label{eqn:1st-term_rhs_Jacobi_1}\nonumber[d_h^\star,\,L_{\omega_{n-1}}]\,d_h\varphi & = & d_h^\star(\omega_{n-1}\wedge(h\partial\varphi^{0,\,1} + \bar\partial\varphi^{1,\,0})) - \omega_{n-1}\wedge d_h^\star d_h\varphi \\
  & = & -i\bar{h}\,d_{-\frac{1}{\bar{h}}}\bigg(h\Lambda_\omega(\partial\varphi^{0,\,1}) + \Lambda_\omega(\bar\partial\varphi^{1,\,0})\bigg)\wedge\omega_{n-1} - d_h^\star d_h\varphi\wedge\omega_{n-1}.\end{eqnarray}

On the other hand, the standard formula (\ref{eqn:prim-form-star-formula-gen}) yields: $$\star\varphi = i(\varphi^{0,\,1} - \varphi^{1,\,0})\wedge\omega_{n-1}.$$ Since $\omega$ is balanced, this implies the first equality on each of the two rows below: \begin{eqnarray*}\partial\star\varphi & = & i\partial(\varphi^{0,\,1} - \varphi^{1,\,0})\wedge\omega_{n-1} = i\partial\varphi^{0,\,1}\wedge\omega_{n-1} = i\Lambda_\omega(\partial\varphi^{0,\,1})\,\omega_n \\
  \bar\partial\star\varphi & = & i\bar\partial(\varphi^{0,\,1} - \varphi^{1,\,0})\wedge\omega_{n-1} = -i\bar\partial\varphi^{1,\,0}\wedge\omega_{n-1} = -i\Lambda_\omega(\bar\partial\varphi^{1,\,0})\,\omega_n.\end{eqnarray*} Taking $-\star$ in each of the above two equalities and using the standard identities $-\star\partial\star = \bar\partial^\star$, $-\star\bar\partial\star = \partial^\star$, we get: \begin{eqnarray}\label{eqn:1st-term_rhs_Jacobi_2}\bar\partial^\star\varphi = -i\Lambda_\omega(\partial\varphi^{0,\,1})  \hspace{3ex} \mbox{and} \hspace{3ex} \partial^\star\varphi = i\Lambda_\omega(\bar\partial\varphi^{1,\,0}).\end{eqnarray}

Putting together (\ref{eqn:1st-term_rhs_Jacobi_1}) and (\ref{eqn:1st-term_rhs_Jacobi_2}), we get: \begin{eqnarray*}[d_h^\star,\,L_{\omega_{n-1}}]\,d_h\varphi & = & -i\bar{h}\,d_{-\frac{1}{\bar{h}}}(ih\bar\partial^\star\varphi - i\partial^\star\varphi)\wedge\omega_{n-1} - d_h^\star d_h\varphi\wedge\omega_{n-1} \\
  & = & h\bar{h}\,d_{-\frac{1}{\bar{h}}}d_{-\frac{1}{\bar{h}}}^\star\varphi\wedge\omega_{n-1} - d_h^\star d_h\varphi\wedge\omega_{n-1}.\end{eqnarray*}

We have thus obtained the following formula: \begin{eqnarray}\label{eqn:1st-term_rhs_Jacobi_conc}[d_h^\star,\,L_{\omega_{n-1}}]\,d_h\varphi & = & \bigg(|h|^2d_{-\frac{1}{\bar{h}}}d_{-\frac{1}{\bar{h}}}^\star - d_h^\star d_h\bigg)\varphi\wedge\omega_{n-1}\end{eqnarray} for every smooth $1$-form $\varphi$ whenever the metric $\omega$ is {\it balanced}.

\vspace{1ex}

$\bullet$ {\it Computation of the second term on the r.h.s. of (\ref{eqn:Delta_h_omega_n-1_Jacobi}) on $1$-forms $\varphi = \varphi^{1,\,0} + \varphi^{0,\,1}$.}

\vspace{1ex}

 We start by computing \begin{eqnarray}\label{eqn:2nd-term_rhs_Jacobi_conc_1}[d_h^\star,\,L_{\omega_{n-1}}]\,\varphi = d_h^\star(\omega_{n-1}\wedge\varphi) - (d_h^\star\varphi)\,\omega_{n-1}.\end{eqnarray}

Since $\omega_{n-1}\wedge\varphi^{1,\,0} = i\star\varphi^{1,\,0}$ and $\omega_{n-1}\wedge\varphi^{0,\,1} = -i\star\varphi^{0,\,1}$, formula (\ref{eqn:d_h-adjoint_star-formula}) yields the first line below: \begin{eqnarray}\label{eqn:2nd-term_rhs_Jacobi_conc_2}\nonumber d_h^\star(\omega_{n-1}\wedge\varphi) & = & i\bar{h}\star d_{\frac{1}{\bar{h}}}(\varphi^{1,\,0}-\varphi^{0,\,1}) = i\star(\partial\varphi^{1,\,0} + \bar{h}\bar\partial\varphi^{1,\,0} - \partial\varphi^{0,\,1} - \bar{h}\bar\partial\varphi^{0,\,1}) \\
  \nonumber & = & \frac{i}{n}\,\Lambda_\omega(\bar{h}\bar\partial\varphi^{1,\,0} - \partial\varphi^{0,\,1})\,\omega_{n-1} -i(\bar{h}\bar\partial\varphi^{1,\,0} - \partial\varphi^{0,\,1})_{prim}\wedge\omega_{n-2} \\
  & + & i(\partial\varphi^{1,\,0} - \bar{h}\bar\partial\varphi^{0,\,1})\wedge\omega_{n-2},\end{eqnarray} where we used the Lefschetz decomposition (\ref{eqn:Lefschetz_decomp_11}) of the $(1,\,1)$-form $\bar{h}\bar\partial\varphi^{1,\,0} - \partial\varphi^{0,\,1}$ and then the standard formula (\ref{eqn:prim-form-star-formula-gen}) to express the value of $\star$ on the primitive forms $\partial\varphi^{1,\,0}$ (of type $(2,\,0)$), $\bar\partial\varphi^{0,\,1}$ (of type $(0,\,2)$) and $(\bar{h}\bar\partial\varphi^{1,\,0} - \partial\varphi^{0,\,1})_{prim}$ (of type $(1,\,1)$) and got $\star(\partial\varphi^{1,\,0}) = \partial\varphi^{1,\,0}\wedge\omega_{n-2}$ and \begin{eqnarray*}\star(\bar\partial\varphi^{0,\,1}) & = & \bar\partial\varphi^{0,\,1}\wedge\omega_{n-2}, \hspace{3ex} \star(\bar{h}\bar\partial\varphi^{1,\,0} - \partial\varphi^{0,\,1})_{prim} = -(\bar{h}\bar\partial\varphi^{1,\,0} - \partial\varphi^{0,\,1})_{prim}\wedge\omega_{n-2}.\end{eqnarray*}

On the other hand, we get \begin{eqnarray}\label{eqn:2nd-term_rhs_Jacobi_conc_3}\nonumber d_h^\star\varphi & = & -\bar{h}\star d_{\frac{1}{\bar{h}}}(\star\varphi^{1,\,0} + \star\varphi^{0,\,1}) = -\bar{h}\star(\frac{1}{\bar{h}}\partial + \bar\partial)(-i\varphi^{1,\,0}\wedge\omega_{n-1} + i\varphi^{0,\,1}\wedge\omega_{n-1}) \\
\nonumber  & \stackrel{(i)}{=} & -\bar{h}\star\bigg(-\frac{i}{\bar{h}}\,\partial\varphi^{1,\,0}\wedge\omega_{n-1} + \frac{i}{\bar{h}}\,\partial\varphi^{0,\,1}\wedge\omega_{n-1} - i\bar\partial\varphi^{1,\,0}\wedge\omega_{n-1} + i\bar\partial\varphi^{0,\,1}\wedge\omega_{n-1}\bigg) \\
\nonumber  & \stackrel{(ii)}{=} &  -\bar{h}\star\bigg(i\bigg(\frac{1}{\bar{h}}\,\partial\varphi^{0,\,1} - \bar\partial\varphi^{1,\,0}\bigg)\wedge\omega_{n-1}\bigg) = -\bar{h}\star\bigg(i\Lambda_\omega\bigg(\frac{1}{\bar{h}}\,\partial\varphi^{0,\,1} - \bar\partial\varphi^{1,\,0}\bigg)\,\omega_n\bigg) \\
  & = & i\Lambda_\omega(\bar{h}\bar\partial\varphi^{1,\,0} - \partial\varphi^{0,\,1}),\end{eqnarray} where the balanced assumption on $\omega$ was used to get (i) and the equalities $\partial\varphi^{1,\,0}\wedge\omega_{n-1} = \bar\partial\varphi^{0,\,1}\wedge\omega_{n-1} = 0$, that hold for bidegree reasons, were used to get (ii).  

Noticing that the last term in (\ref{eqn:2nd-term_rhs_Jacobi_conc_3}) also features within the first term on the second line in (\ref{eqn:2nd-term_rhs_Jacobi_conc_2}), the conclusion of (\ref{eqn:2nd-term_rhs_Jacobi_conc_2}) can be re-written as \begin{eqnarray}\label{eqn:2nd-term_rhs_Jacobi_conc_3_after}d_h^\star(\omega_{n-1}\wedge\varphi) = \frac{1}{n}\,(d_h^\star\varphi)\,\omega_{n-1} -i(\bar{h}\bar\partial\varphi^{1,\,0} - \partial\varphi^{0,\,1})_{prim}\wedge\omega_{n-2} + i(\partial\varphi^{1,\,0} - \bar{h}\bar\partial\varphi^{0,\,1})\wedge\omega_{n-2}.\end{eqnarray} From this and from (\ref{eqn:2nd-term_rhs_Jacobi_conc_1}), we get: \begin{eqnarray*}[d_h^\star,\,L_{\omega_{n-1}}]\,\varphi = \bigg(\frac{1}{n}-1\bigg)\,(d_h^\star\varphi)\,\omega_{n-1} -i(\bar{h}\bar\partial\varphi^{1,\,0} - \partial\varphi^{0,\,1})_{prim}\wedge\omega_{n-2} + i(\partial\varphi^{1,\,0} - \bar{h}\bar\partial\varphi^{0,\,1})\wedge\omega_{n-2}.\end{eqnarray*} Hence, using the balanced hypothesis $d_h\omega=0$, we get the first two lines below: \begin{eqnarray*}d_h[d_h^\star,\,L_{\omega_{n-1}}]\,\varphi & = & \bigg(\frac{1}{n}-1\bigg)\,dd_h^\star\varphi\wedge\omega_{n-1} -id_h\bigg((\bar{h}\bar\partial\varphi^{1,\,0} - \partial\varphi^{0,\,1})\wedge\omega_{n-2}\bigg) \\
  & + & \frac{n-1}{n}i\,d_h\bigg(\Lambda_\omega(\bar{h}\bar\partial\varphi^{1,\,0} - \partial\varphi^{0,\,1})\bigg)\wedge\omega_{n-1} \\
  & + & i(\partial\varphi^{1,\,0} - \bar{h}\bar\partial\varphi^{0,\,1})\wedge d_h\omega_{n-2} - i(|h|^2\,\partial\bar\partial\varphi^{0,\,1} + \partial\bar\partial\varphi^{1,\,0})\wedge\omega_{n-2}.\end{eqnarray*} Now, formula (\ref{eqn:2nd-term_rhs_Jacobi_conc_3}) shows that the term on the second line above equals minus the first term on the r.h.s. of the first line. Hence, the sum of these two terms vanishes and we get: \begin{eqnarray}\label{eqn:2nd-term_rhs_Jacobi_conc}\nonumber d_h[d_h^\star,\,L_{\omega_{n-1}}]\,\varphi & = & i\,\bigg(\partial\varphi^{1,\,0} + (\partial\varphi^{0,\,1} - \bar{h}\bar\partial\varphi^{1,\,0}) - \bar{h}\bar\partial\varphi^{0,\,1}\bigg)\wedge d_h\omega_{n-2} -i(|h|^2 + 1)\,\partial\bar\partial\varphi\wedge\omega_{n-2} \\
 & = & -i\bar{h}\,d_{-\frac{1}{\bar{h}}}\varphi\wedge d_h\omega_{n-2} -i(|h|^2 + 1)\,\partial\bar\partial\varphi\wedge\omega_{n-2}.\end{eqnarray}

\vspace{1ex}

$\bullet$ {\it Conclusion.}

\vspace{1ex}

Putting together (\ref{eqn:Delta_h_omega_n-1_Jacobi}), (\ref{eqn:1st-term_rhs_Jacobi_conc}) and (\ref{eqn:2nd-term_rhs_Jacobi_conc}), we get (\ref{eqn:commutation-defect_1-forms_pre-integration}). The proof of Lemma \ref{Lem:commutation-defect_1-forms_pre-integration} is complete.  \hfill $\Box$

\vspace{3ex}

Recall that for any Hermitian metric $\omega$ on an $n$-dimensional complex manifold $X$, the pointwise {\it Lefschetz map}: \begin{equation*}L_{\omega_{n-1}}:\Lambda^1T^\star X\longrightarrow\Lambda^{2n-1}T^\star X, \hspace{3ex} \varphi\longmapsto\psi:=\omega_{n-1}\wedge\varphi,\end{equation*} is {\it bijective} and a {\it quasi-isometry} (in the sense of Lemma \ref{Lem:appendix_pointwise_quasi-isometry}).

We will now integrate the result of Lemma \ref{Lem:commutation-defect_1-forms_pre-integration} expressing the commutation defect between $\Delta_h$ and $L_{\omega_{n-1}}$ on $1$-forms. We need to assume our balanced metric $\omega$ to be {\it complete} to ensure that the two meanings of $d_h^\star$ coincide and the $L^2_\omega$-inner products can be handled as in the compact case (see (b) and (c) of Theorem \ref{The:standard_complete_operators}).

\begin{Prop}\label{Prop:commutation-defect_1-forms_integrated} Let $X$ be a complex manifold with $\mbox{dim}_\C X=n$. Suppose there exists a {\bf complete balanced metric} $\omega$ on $X$.

  Then, for any $h\in\C^\star$ and any {\bf $1$-form} $\varphi\in\mbox{Dom}\,(\Delta_{-\frac{1}{\bar{h}}})$ on $X$, the following identity holds: \begin{eqnarray}\label{eqn:commutation-defect_1-forms_integrated}\langle\langle\Delta_h(\omega_{n-1}\wedge\varphi),\,\omega_{n-1}\wedge\varphi\rangle\rangle = |h|^2\, \langle\langle\Delta_{-\frac{1}{\bar{h}}}\varphi,\,\varphi\rangle\rangle.\end{eqnarray}

\end{Prop}

\noindent {\it Proof.} Throughout the proof, $\varphi$ will stand for an arbitrary smooth $1$-form on $X$.

\vspace{1ex}

$\bullet$ We first notice that $d_hd_{-\frac{1}{\bar{h}}}\varphi = ((|h|^2+1)/\bar{h})\,\partial\bar\partial\varphi$, hence \begin{eqnarray*}d_h\bigg(-i\bar{h}\,d_{-\frac{1}{\bar{h}}}\varphi\wedge\omega_{n-2}\bigg) = -i\bar{h}\,d_{-\frac{1}{\bar{h}}}\varphi\wedge d_h\omega_{n-2} -i(|h|^2 + 1)\,\partial\bar\partial\varphi\wedge\omega_{n-2}.\end{eqnarray*} These are the last two terms of formula (\ref{eqn:commutation-defect_1-forms_pre-integration}).

Putting $\psi:=\omega_{n-1}\wedge\varphi$ and using (\ref{eqn:commutation-defect_1-forms_pre-integration}) with its last two terms transformed as above, we get: \begin{eqnarray}\label{eqn:commutation-defect_1-forms_integrated_1}\nonumber\langle\langle\Delta_h\psi,\,\psi\rangle\rangle & = & \langle\langle\Delta_h\varphi\wedge\omega_{n-1},\,\varphi\wedge\omega_{n-1}\rangle\rangle + \langle\langle(|h|^2d_{-\frac{1}{\bar{h}}}d_{-\frac{1}{\bar{h}}}^\star - d_h^\star d_h)\,\varphi\wedge\omega_{n-1},\,\varphi\wedge\omega_{n-1}\rangle\rangle \\
  \nonumber &  & \hspace{27ex} - i\bar{h}\,\langle\langle d_{-\frac{1}{\bar{h}}}\varphi\wedge\omega_{n-2},\,d_h^\star(\omega_{n-1}\wedge\varphi)\rangle\rangle  \\
\nonumber  & = & \langle\langle d_hd_h^\star\varphi\wedge\omega_{n-1},\,\varphi\wedge\omega_{n-1}\rangle\rangle + \langle\langle |h|^2d_{-\frac{1}{\bar{h}}}d_{-\frac{1}{\bar{h}}}^\star\varphi\wedge\omega_{n-1},\,\varphi\wedge\omega_{n-1}\rangle\rangle \\
\nonumber &  & \hspace{27ex} - i\bar{h}\,\langle\langle d_{-\frac{1}{\bar{h}}}\varphi\wedge\omega_{n-2},\,d_h^\star(\omega_{n-1}\wedge\varphi)\rangle\rangle  \\
\nonumber & \stackrel{(i)}{=} & \langle\langle d_hd_h^\star\varphi,\,\varphi\rangle\rangle + |h|^2\,\langle\langle d_{-\frac{1}{\bar{h}}}d_{-\frac{1}{\bar{h}}}^\star\varphi,\,\varphi\rangle\rangle - i\bar{h}\,\langle\langle d_{-\frac{1}{\bar{h}}}\varphi\wedge\omega_{n-2},\,d_h^\star(\omega_{n-1}\wedge\varphi)\rangle\rangle  \\
  & = & ||d_h^\star\varphi||^2 + |h|^2\,||d_{-\frac{1}{\bar{h}}}^\star\varphi||^2 - i\bar{h}\,\langle\langle d_{-\frac{1}{\bar{h}}}\varphi\wedge\omega_{n-2},\,d_h^\star(\omega_{n-1}\wedge\varphi)\rangle\rangle,\end{eqnarray} where (i) followed from Lemma \ref{Lem:appendix_pointwise-isometry_primitive} applied to the (necessarily primitive) $1$-forms $\varphi$, $d_hd_h^\star\varphi$ and $d_{-\frac{1}{\bar{h}}}d_{-\frac{1}{\bar{h}}}^\star\varphi$.

\vspace{1ex}

$\bullet$ We now transform the last term in (\ref{eqn:commutation-defect_1-forms_integrated_1}), namely $T(\varphi):=- i\bar{h}\,\langle\langle d_{-\frac{1}{\bar{h}}}\varphi\wedge\omega_{n-2},\,d_h^\star(\omega_{n-1}\wedge\varphi)\rangle\rangle$.

\vspace{1ex}

Since the multiplication map $\omega_{n-2}\wedge\cdot : \Lambda^2T^\star X\longrightarrow\Lambda^{2n-2}T^\star X$ is {\it bijective}, there exists a unique $2$-form $\beta$ such that $d_h^\star(\omega_{n-1}\wedge\varphi) = \omega_{n-2}\wedge\beta$. Thus, using (\ref{eqn:2nd-term_rhs_Jacobi_conc_3_after}) for the second equality below, we get: \begin{eqnarray*}\omega_{n-2}\wedge\beta = d_h^\star(\omega_{n-1}\wedge\varphi) = \omega_{n-2}\wedge\bigg(\frac{1}{n(n-1)}\,(d_h^\star\varphi)\,\omega -i(\bar{h}\bar\partial\varphi^{1,\,0} - \partial\varphi^{0,\,1})_{prim} + i(\partial\varphi^{1,\,0} - \bar{h}\bar\partial\varphi^{0,\,1})\bigg).\end{eqnarray*} The uniqueness of $\beta$ implies that \begin{eqnarray}\label{eqn:commutation-defect_1-forms_integrated_2}\beta = \bigg(-i(\bar{h}\bar\partial\varphi^{1,\,0} - \partial\varphi^{0,\,1})_{prim} + i(\partial\varphi^{1,\,0} - \bar{h}\bar\partial\varphi^{0,\,1})\bigg) + \frac{1}{n(n-1)}\,(d_h^\star\varphi)\,\omega.\end{eqnarray} In particular, the primitive part $\beta_{prim}$ of $\beta$ in the Lefschetz decomposition is the form inside the large parenthesis and $\Lambda_\omega\beta = \frac{1}{n-1}\,d_h^\star\varphi$.

On the other hand, we have \begin{eqnarray}\label{eqn:commutation-defect_1-forms_integrated_3}d_{-\frac{1}{\bar{h}}}\varphi = -\frac{1}{\bar{h}}\,\partial\varphi^{1,\,0} + (-\frac{1}{\bar{h}}\,\partial\varphi^{0,\,1} + \bar\partial\varphi^{1,\,0})_{prim} + \bar\partial\varphi^{0,\,1} + \frac{1}{ni\bar{h}}\,(d_h^\star\varphi)\,\omega,\end{eqnarray} where the value of the last term follows from formula (\ref{eqn:2nd-term_rhs_Jacobi_conc_3}). This implies that $\beta$ and $-i\bar{h}\,d_{-\frac{1}{\bar{h}}}\varphi$ have the same primitive part: \begin{eqnarray}\label{eqn:commutation-defect_1-forms_integrated_4}\beta_{prim} = -i\bar{h}\,(d_{-\frac{1}{\bar{h}}}\varphi)_{prim}.\end{eqnarray}

We get: \begin{eqnarray*}\langle d_{-\frac{1}{\bar{h}}}\varphi\wedge\omega_{n-2},\,d_h^\star(\omega_{n-1}\wedge\varphi)\rangle & = & \langle d_{-\frac{1}{\bar{h}}}\varphi\wedge\omega_{n-2},\,\beta\wedge\omega_{n-2}\rangle \\
  & = & \langle(d_{-\frac{1}{\bar{h}}}\varphi)_{prim},\,\beta_{prim}\rangle + (n-1)^2n\,\bigg\langle\frac{1}{ni\bar{h}}\,d_h^\star\varphi,\,\frac{1}{n(n-1)}\,d_h^\star\varphi\bigg\rangle,\end{eqnarray*} where the last equality follows from formula (\ref{eqn:appendix_pointwise_quasi-isometry_1}) in the Appendix.

From this and from (\ref{eqn:commutation-defect_1-forms_integrated_2})-(\ref{eqn:commutation-defect_1-forms_integrated_4}), we get: \begin{eqnarray}\label{eqn:commutation-defect_1-forms_integrated_5}\nonumber T(\varphi) & = & - i\bar{h}\,\langle\langle d_{-\frac{1}{\bar{h}}}\varphi\wedge\omega_{n-2},\,d_h^\star(\omega_{n-1}\wedge\varphi)\rangle\rangle \\
  & = & ||\partial\varphi^{1,\,0}||^2 + ||(\partial\varphi^{0,\,1} - \bar{h}\,\bar\partial\varphi^{1,\,0})_{prim}||^2 + |h|^2\, ||\bar\partial\varphi^{0,\,1}||^2 -(1-\frac{1}{n})\,||d_h^\star\varphi||^2.\end{eqnarray}

\vspace{1ex}

$\bullet$ Putting (\ref{eqn:commutation-defect_1-forms_integrated_1}) and (\ref{eqn:commutation-defect_1-forms_integrated_5}) together and writing $\frac{1}{n}\,||d_h^\star\varphi||^2 = |h|^2\,\frac{1}{n}\,||\frac{1}{i\bar{h}}\,d_h^\star\varphi||^2$, we get: \begin{eqnarray*}\langle\langle\Delta_h\psi,\,\psi\rangle\rangle & = & |h|^2\,||d_{-\frac{1}{\bar{h}}}^\star\varphi||^2 + |h|^2\,\frac{1}{n}\,||\frac{1}{i\bar{h}}\,d_h^\star\varphi||^2  \\
  & + & |h|^2\,\bigg(||-\frac{1}{\bar{h}}\,\partial\varphi^{1,\,0}||^2 + ||(-\frac{1}{\bar{h}}\,\partial\varphi^{0,\,1} + \bar\partial\varphi^{1,\,0})_{prim}||^2 + ||\bar\partial\varphi^{0,\,1}||^2\bigg).\end{eqnarray*} Thanks to the expression of $d_{-\frac{1}{\bar{h}}}\varphi$ obtained in (\ref{eqn:commutation-defect_1-forms_integrated_3}), this translates to \begin{eqnarray*}\langle\langle\Delta_h\psi,\,\psi\rangle\rangle = |h|^2\,\bigg(||d_{-\frac{1}{\bar{h}}}^\star\varphi||^2 + ||d_{-\frac{1}{\bar{h}}}\varphi||^2\bigg) = |h|^2\,\langle\langle\Delta_{-\frac{1}{\bar{h}}}\varphi,\,\varphi\rangle\rangle,\end{eqnarray*} which is (\ref{eqn:commutation-defect_1-forms_integrated}).

Proposition \ref{Prop:commutation-defect_1-forms_integrated} is proved.  \hfill $\Box$

\vspace{3ex}

An immediate consequence of Proposition \ref{Prop:commutation-defect_1-forms_integrated} is the following {\it Hard Lefschetz}-type result for spaces of harmonic $L^2_\omega$-forms induced by a given {\it complete balanced metric} $\omega$ and different operators $\Delta_{-\frac{1}{\bar{h}}}$ and $\Delta_h$. Note that $h\neq -\frac{1}{\bar{h}}$ for all $h\in\C^\star$. This is the price we have to pay in the non-K\"ahler balanced context to get this kind of results.

\begin{Cor}\label{Cor:bal-complete_Hard-Lefschetz} Let $X$ be a complex manifold with $\mbox{dim}_\C X=n$. Suppose there exists a {\bf complete balanced metric} $\omega$ on $X$. Then, for any $h\in\C^\star$, the map \begin{equation*}\omega_{n-1}\wedge\cdot:{\cal H}^1_{\Delta_{-\frac{1}{\bar{h}}}}(X,\,\C)\longrightarrow{\cal H}^{2n-1}_{\Delta_h}(X,\,\C),  \hspace{3ex} \varphi\longmapsto\omega_{n-1}\wedge\varphi,\end{equation*} is {\bf well-defined} and an {\bf isomorphism}. 

\end{Cor}  

\noindent {\it Proof.} The well-definedness, namely the fact that this map takes $\Delta_{-\frac{1}{\bar{h}}}$-harmonic $L^2_\omega$-forms to $\Delta_h$-harmonic $L^2_\omega$-forms, follows at once from Proposition \ref{Prop:commutation-defect_1-forms_integrated} and from the form $\omega_{n-1}$ being $\omega$-bounded. The fact that this map is an isomorphism follows from the standard fact that the corresponding pointwise map is bijective.  \hfill $\Box$

\begin{Cor}\label{Cor:bal-complete-exact_Laplacian_l-bound} Let $X$ be a complex manifold with $\mbox{dim}_\C X=n$. Suppose there exists a {\bf complete balanced metric} $\omega$ on $X$ such that $\omega_{n-1} = d\Gamma$ for an {\bf $\omega$-bounded} smooth $(2n-3)$-form $\Gamma$. Then \begin{equation}\label{eqn:bal-complete-exact_Laplacian_l-bound}\langle\langle\Delta\psi,\,\psi\rangle\rangle\geq\frac{1}{4||\Gamma||^2_{L^\infty_\omega}}\,||\psi||^2\end{equation} for every {\bf pure-type} form $\psi\in\mbox{Dom}(\Delta)$ of degree $2n-1$.

\end{Cor} 

\noindent {\it Proof.} Taking $h=1$ in Proposition \ref{Prop:commutation-defect_1-forms_integrated}, (\ref{eqn:commutation-defect_1-forms_integrated}) gives: \begin{eqnarray*}\langle\langle\Delta\psi,\,\psi\rangle\rangle = \langle\langle\Delta_{-1}\varphi,\,\varphi\rangle\rangle = ||(\partial-\bar\partial)\varphi ||^2 + ||(\partial-\bar\partial)^\star\varphi ||^2\geq||(\partial-\bar\partial)\varphi ||^2,\end{eqnarray*} for every $(2n-1)$-form $\psi$, where $\varphi$ is the unique $1$-form such that $\psi=\omega_{n-1}\wedge\varphi$. (See isomorphism (\ref{eqn:pointwise-Lefschetz-map}) for $r=1$.) Meanwhile, $\psi$ is of pure type (either $(n,\,n-1)$ or $(n-1,\,n)$) if and only if $\varphi$ is of pure type (respectively, either $(1,\,0)$ or $(0,\,1)$). In this case, $\partial\varphi$ and $\bar\partial\varphi$ are of different pure types, hence orthogonal to each other, hence $||(\partial-\bar\partial)\varphi ||^2 = ||(\partial+\bar\partial)\varphi ||^2$. Thus, we get: \begin{eqnarray}\label{eqn:bal-complete-exact_Laplacian_l-bound_1}\langle\langle\Delta\psi,\,\psi\rangle\rangle\geq||d\varphi ||^2,\end{eqnarray} for every {\bf pure-type} $(2n-1)$-form $\psi\in\mbox{Dom}(\Delta)$.

To complete the proof, we adapt the proof of Theorem 1.4.A. in [Gro91] to our context.

\vspace{1ex}

 Since any $1$-form $\varphi$ is primitive, Lemma \ref{Lem:appendix_pointwise-isometry_primitive} gives: $|\psi|^2 = |\omega_{n-1}\wedge\varphi|^2 = |\varphi|^2$. In particular, \begin{eqnarray}\label{eqn:bal-complete-exact_Laplacian_l-bound_2}||\psi|| = ||\varphi||.\end{eqnarray}

 Meanwhile, we have: $\psi= \omega_{n-1}\wedge\varphi = d\Gamma\wedge\varphi  = d(\Gamma\wedge\varphi) + \Gamma\wedge d\varphi$. In other words, \begin{eqnarray}\label{eqn:bal-complete-exact_Laplacian_l-bound_3}\psi= d\theta + \psi', \hspace{5ex} \mbox{where}\hspace{1ex} \theta:=\Gamma\wedge\varphi \hspace{2ex}\mbox{and}\hspace{2ex} \psi':=\Gamma\wedge d\varphi.\end{eqnarray}

 To estimate $\theta$, we write: \begin{eqnarray}\label{eqn:bal-complete-exact_Laplacian_l-bound_4}||\theta||\leq||\Gamma||_{L^\infty_\omega}\,||\varphi|| = ||\Gamma||_{L^\infty_\omega}\,||\psi||,\end{eqnarray} where (\ref{eqn:bal-complete-exact_Laplacian_l-bound_2}) was used to get the last equality.

 To estimate $\psi'$, we write: \begin{eqnarray}\label{eqn:bal-complete-exact_Laplacian_l-bound_5}||\psi'||\leq||\Gamma||_{L^\infty_\omega}\,||d\varphi|| \leq ||\Gamma||_{L^\infty_\omega}\,\langle\langle\Delta\psi,\,\psi\rangle\rangle^{\frac{1}{2}},\end{eqnarray} where (\ref{eqn:bal-complete-exact_Laplacian_l-bound_1}) and the fact that $\varphi$ is of pure type were used to get the last inequality.
 
 To find an upper bound for $||\psi||$, we write: \begin{eqnarray}\label{eqn:bal-complete-exact_Laplacian_l-bound_6}||\psi||^2 = \langle\langle\psi,\,d\theta + \psi'\rangle\rangle\leq|\langle\langle\psi,\,d\theta\rangle\rangle| + |\langle\langle\psi,\,\psi'\rangle\rangle|,\end{eqnarray} where (\ref{eqn:bal-complete-exact_Laplacian_l-bound_3}) was used to get the first equality.

 For the first term on the r.h.s. of (\ref{eqn:bal-complete-exact_Laplacian_l-bound_6}), we get: \begin{eqnarray}\label{eqn:bal-complete-exact_Laplacian_l-bound_7}|\langle\langle\psi,\,d\theta\rangle\rangle|=|\langle\langle d^\star\psi,\,\theta\rangle\rangle|\leq||d^\star\psi||\,||\theta||\leq\langle\langle\Delta\psi,\,\psi\rangle\rangle^{\frac{1}{2}}\,||\Gamma||_{L^\infty_\omega}\,||\psi||,\end{eqnarray} where (\ref{eqn:bal-complete-exact_Laplacian_l-bound_4}) was used to get the last inequality.

 For the second term on the r.h.s. of (\ref{eqn:bal-complete-exact_Laplacian_l-bound_6}), we get: \begin{eqnarray}\label{eqn:bal-complete-exact_Laplacian_l-bound_8}|\langle\langle\psi,\,\psi'\rangle\rangle|\leq||\psi'||\,||\psi||\leq||\Gamma||_{L^\infty_\omega}\,\langle\langle\Delta\psi,\,\psi\rangle\rangle^{\frac{1}{2}}\,||\psi||,\end{eqnarray}  where (\ref{eqn:bal-complete-exact_Laplacian_l-bound_5}) was used to get the last inequality.

 Adding up (\ref{eqn:bal-complete-exact_Laplacian_l-bound_7}) and (\ref{eqn:bal-complete-exact_Laplacian_l-bound_8}) and using (\ref{eqn:bal-complete-exact_Laplacian_l-bound_6}), we get \begin{eqnarray*}||\psi||\leq 2\,||\Gamma||_{L^\infty_\omega}\,\langle\langle\Delta\psi,\,\psi\rangle\rangle^{\frac{1}{2}},\end{eqnarray*} which is (\ref{eqn:bal-complete-exact_Laplacian_l-bound}). The proof is complete. \hfill $\Box$

 \vspace{3ex}

 For the record, if we do not assume $\psi$ to be of pure type and use the full force of (\ref{eqn:commutation-defect_1-forms_integrated}) rather than (\ref{eqn:bal-complete-exact_Laplacian_l-bound_1}), we can run the argument in the proof of Corollary \ref{Cor:bal-complete-exact_Laplacian_l-bound} with minor modifications starting from the observation that $\omega_{n-1} = d_{-\frac{1}{\bar{h}}}\Gamma_{-\frac{1}{\bar{h}}}$, where $\Gamma_h:=h\,\Gamma^{n,\,n-3} + \Gamma^{n-1,\,n-2} + (1/h)\,\Gamma^{n-2,\,n-1} + (1/h^2)\,\Gamma^{n-3,\,n}$ for every $h\in\C^\star$ and the $\Gamma^{p,\,q}$'s are the pure-type components of $\Gamma$. Then, we get the following analogue of (\ref{eqn:bal-complete-exact_Laplacian_l-bound}): \begin{equation}\label{eqn:bal-complete-exact_Laplacian_l-bound_h}||\psi||\leq C_h\,||\Gamma_{-\frac{1}{\bar{h}}}||\,\bigg(\langle\langle\Delta_h\psi,\,\psi\rangle\rangle^{\frac{1}{2}} + \langle\langle\Delta_{-\frac{1}{\bar{h}}}\psi,\,\psi\rangle\rangle^{\frac{1}{2}}\bigg)\end{equation} for every form $\psi\in\mbox{Dom}(\Delta_h)\cap\mbox{Dom}(\Delta_{-\frac{1}{\bar{h}}})$ (not necessarily of pure type) of degree $2n-1$, where $C_h:=\max(1, 1/|h|)$.

 The occurrence of two different Laplacians on the r.h.s. of (\ref{eqn:bal-complete-exact_Laplacian_l-bound_h}) (recall that $h\neq -\frac{1}{\bar{h}}$ for every $h\in\C^\star)$ is the downside of that estimate that we avoided in Corollary \ref{Cor:bal-complete-exact_Laplacian_l-bound} by restricting attention to pure-type forms. The advantage of dealing with a single Laplacian is demonstrated by Theorem \ref{The:bal-complete-exact_no-harmonic-forms_degree1} in the introduction that we now prove as a consequence of the above discussion.

%\begin{Cor}\label{Cor:bal-complete-exact_no-harmonic-forms_degree1} Let $X$ be a compact complex {\bf balanced hyperbolic} manifold with $\mbox{dim}_\C X=n$. Let $\pi:\widetilde{X}\longrightarrow X$ be the universal cover of $X$ and $\widetilde\omega:=\pi^\star\omega$ the lift to $\widetilde{X}$ of a balanced hyperbolic metric $\omega$ on $X$.

%   There are no non-zero $\Delta_{\widetilde\omega}$-harmonic $L^2_{\widetilde\omega}$-forms of pure types and of degrees $1$ and $2n-1$ on $\widetilde{X}$: $${\cal H}^{1,\,0}_{\Delta_{\widetilde\omega}}(\widetilde{X},\,\C) = {\cal H}^{0,\,1}_{\Delta_{\widetilde\omega}}(\widetilde{X},\,\C) = 0  \hspace{3ex} \mbox{and} \hspace{3ex} {\cal H}^{n,\,n-1}_{\Delta_{\widetilde\omega}}(\widetilde{X},\,\C) = {\cal H}^{n-1,\,n}_{\Delta_{\widetilde\omega}}(\widetilde{X},\,\C) = 0.$$

 %\end{Cor}

\vspace{3ex} 

\noindent {\it Proof of Theorem \ref{The:bal-complete-exact_no-harmonic-forms_degree1}.} The pair $(\widetilde{X},\, \widetilde\omega)$ satisfies the hypotheses of Corollary \ref{Cor:bal-complete-exact_Laplacian_l-bound} (playing the role of the pair $(X,\,\omega)$ therein). When applied to $(n,\,n-1)$-forms and to $(n-1,\,n)$-forms $\psi\in\mbox{Dom}(\Delta_{\widetilde\omega})$, inequality (\ref{eqn:bal-complete-exact_Laplacian_l-bound}) gives the following implication: $$\Delta_{\widetilde\omega}\psi=0 \implies \psi=0.$$
This proves the vanishing of ${\cal H}^{n,\,n-1}_{\Delta_{\widetilde\omega}}(\widetilde{X},\,\C)$ and ${\cal H}^{n-1,\,n}_{\Delta_{\widetilde\omega}}(\widetilde{X},\,\C)$.

Meanwhile, the Hodge star operator $\star = \star_{\widetilde\omega}$ commutes with $\Delta_{\widetilde\omega}$, so it induces isomorphisms $$\star_{\widetilde\omega}:{\cal H}^{1,\,0}_{\Delta_{\widetilde\omega}}(\widetilde{X},\,\C)\longrightarrow{\cal H}^{n,\,n-1}_{\Delta_{\widetilde\omega}}(\widetilde{X},\,\C) \hspace{3ex} \mbox{and} \hspace{3ex} \star_{\widetilde\omega}:{\cal H}^{0,\,1}_{\Delta_{\widetilde\omega}}(\widetilde{X},\,\C)\longrightarrow{\cal H}^{n-1,\,n}_{\Delta_{\widetilde\omega}}(\widetilde{X},\,\C).$$ Therefore, the spaces ${\cal H}^{1,\,0}_{\Delta_{\widetilde\omega}}(\widetilde{X},\,\C)$ and ${\cal H}^{0,\,1}_{\Delta_{\widetilde\omega}}(\widetilde{X},\,\C)$ must vanish as well. \hfill $\Box$

\subsection{Harmonic $L^2$-forms of degree $2$ on the universal cover of a balanced hyperbolic manifold}\label{subsection:harmonic_bal-hyp_degree2}

We will discuss $2$-forms in a way analogous to the discussion of $1$-forms we had in $\S.$\ref{subsection:harmonic_bal-hyp_degree1}. The context and the notation are the same. The analogue of Lemma \ref{Lem:commutation-defect_1-forms_pre-integration} is

\begin{Lem}\label{Lem:commutation-defect_2-forms_pre-integration} Let $X$ be a complex manifold with $\mbox{dim}_\C X=n$. Suppose there exists a {\bf balanced metric} $\omega$ on $X$. Then, for any $h\in\C^\star$ and any {\bf $2$-form} $\alpha$ on $X$, the following identity holds: \begin{eqnarray}\label{eqn:commutation-defect_2-forms_pre-integration}[\Delta_h,\,L_{\omega_{n-1}}]\alpha = -(|h|^2 + 1)\,i\partial\bar\partial(\Lambda_\omega\alpha)\wedge\omega_{n-1} - \omega_{n-1}\wedge\Delta_h\alpha.\end{eqnarray}

\end{Lem}  

\noindent {\it Proof.} We compute separately the two terms applied to $\alpha$ on the r.h.s. of the consequence (\ref{eqn:Delta_h_omega_n-1_Jacobi}) of the Jacobi identity and the balanced hypothesis on $\omega$.

The first term is \begin{eqnarray}\label{eqn:commutation-defect_2-forms_pre-integration_1}[d_h^\star,\,L_{\omega_{n-1}}]\,d_h\alpha = - \omega_{n-1}\wedge d_h^\star d_h\alpha,\end{eqnarray} since $d_h^\star(\omega_{n-1}\wedge d_h\alpha)=0$ owing to the vanishing of $\omega_{n-1}\wedge d_h\alpha$ for degree reasons.

To compute $d_h\,[d_h^\star,\,L_{\omega_{n-1}}]\alpha$, we notice that \begin{eqnarray*}[d_h^\star,\,L_{\omega_{n-1}}]\alpha = d_h^\star(\omega_{n-1}\wedge\alpha) - \omega_{n-1}\wedge d_h^\star\alpha = -i\bar{h}\,d_{-\frac{1}{\bar{h}}}(\Lambda_\omega\alpha)\wedge\omega_{n-1} - \omega_{n-1}\wedge d_h^\star\alpha,\end{eqnarray*} where the last identity follows from (\ref{eqn:d_h_star_omegan-1_11-form}). Thus, using the balanced hypothesis on $\omega$, we get: \begin{eqnarray}\label{eqn:commutation-defect_2-forms_pre-integration_2}d_h\,[d_h^\star,\,L_{\omega_{n-1}}]\alpha = -i\bar{h}d_h\,d_{-\frac{1}{\bar{h}}}(\Lambda_\omega\alpha)\wedge\omega_{n-1} - \omega_{n-1}\wedge d_hd_h^\star\alpha.\end{eqnarray}

Finally, $d_h\,d_{-\frac{1}{\bar{h}}} = ((|h|^2 + 1)/\bar{h})\,\partial\bar\partial$, so (\ref{eqn:commutation-defect_2-forms_pre-integration}) follows from (\ref{eqn:commutation-defect_2-forms_pre-integration_1}) and (\ref{eqn:commutation-defect_2-forms_pre-integration_2}).  \hfill $\Box$

\vspace{3ex}

We now deduce the following analogue of Proposition \ref{Prop:commutation-defect_1-forms_integrated}.

\begin{Prop}\label{Prop:commutation-defect_2-forms_integrated} Let $(X,\,\omega)$ be a {\bf complete balanced} manifold, $\mbox{dim}_{\C}X=n\geq 2$.

For any $h\in\C^\star$ and any {\bf $2$-form} $\varphi\in\mbox{Dom}\,(\Delta_h)$ on $X$, the following identity holds: \begin{eqnarray}\label{eqn:commutation-defect_2-forms_integrated}\langle\langle\Delta_h(\omega_{n-1}\wedge\alpha),\,\omega_{n-1}\wedge\alpha\rangle\rangle = (|h|^2 + 1)\,||\bar\partial(\Lambda_\omega\alpha)||^2 .\end{eqnarray}

\end{Prop}

\noindent {\it Proof.} An immediate consequence of (\ref{eqn:commutation-defect_2-forms_pre-integration}) is the identity $$\Delta_h(\omega_{n-1}\wedge\alpha) = -(|h|^2 + 1)\,i\partial\bar\partial(\Lambda_\omega\alpha)\wedge\omega_{n-1}.$$ Taking the pointwise inner product (w.r.t. $\omega$) against $\omega_{n-1}\wedge\alpha$ and using the Lefschetz decomposition $\alpha^{1,\,1} = \alpha^{1,\,1}_{prim} + (1/n)\,(\Lambda_\omega\alpha^{1,\,1})\,\omega$ of the $(1,\,1)$-type component of $\alpha$, its analogue for the $(1,\,1)$-form $i\partial\bar\partial(\Lambda_\omega\alpha)$ and the fact that the product of any primitive $2$-form with $\omega_{n-1}$ vanishes, we get: \begin{eqnarray}\label{eqn:commutation-defect_2-forms_integrated_1}\nonumber\langle\Delta_h(\omega_{n-1}\wedge\alpha),\,\omega_{n-1}\wedge\alpha\rangle & = & -(|h|^2 + 1)\,\langle\widetilde\Delta_\omega(\Lambda_\omega\alpha)\,\omega_n,\,(\Lambda_\omega\alpha)\,\omega_n\rangle \\
   & = & -(|h|^2 + 1)\,\langle\widetilde\Delta_\omega(\Lambda_\omega\alpha),\,\Lambda_\omega\alpha\rangle,\end{eqnarray} where $\widetilde\Delta_\omega f:= \Lambda_\omega(i\partial\bar\partial f)$ for any function $f$ on $X$. It is standard that the Laplacian $\widetilde\Delta_\omega$ is a non-positive operator on functions.
 Identity (\ref{eqn:appendix_pointwise-isometry_primitive}) in Lemma \ref{Lem:appendix_pointwise-isometry_primitive} with $k=0$ and $r=n$ was used to get the last equality in (\ref{eqn:commutation-defect_2-forms_integrated_1}). 

 Now, we need the following simple observation.

 \begin{Lem}\label{Lem:Laplacian-balanced_inner-rpod} Let $(X,\,\omega)$ be a {\bf complete balanced} manifold, $\mbox{dim}_{\C}X=n\geq 2$. For any function $f\in\mbox{Dom}\,(\widetilde\Delta_\omega)$, we have: $\langle\langle\widetilde\Delta_\omega f,\,f\rangle\rangle = -||\bar\partial f||^2$.

\end{Lem}

\noindent {\it Proof of Lemma \ref{Lem:Laplacian-balanced_inner-rpod}.} The formula $\partial^{\star} = -\star\bar\partial\star$ gives the third equality below: \begin{eqnarray*}\label{eqn:delta-phi1}\nonumber\langle\langle\widetilde\Delta_{\omega}f,\,f\rangle\rangle & = & \langle\langle\Lambda_{\omega}(i\partial\bar\partial f),\,f\rangle\rangle = \langle\langle i\bar\partial f,\,\partial^{\star}(f\omega)\rangle\rangle = -i\,\langle\langle\bar\partial f,\,\star\,\bar\partial(f\,\omega_{n-1})\rangle\rangle \\
   & = & - i\,\langle\langle\bar\partial f,\, \star(\bar\partial f\wedge\omega_{n-1})\rangle\rangle,\end{eqnarray*} where we used the balanced hypothesis on $\omega$ to get the last equality.

 Now, $\bar\partial f$ is a $(0,\, 1)$-form, hence {\it primitive}, so the standard formula (\ref{eqn:prim-form-star-formula-gen}) yields: $$\star(i\bar\partial f) = -\bar\partial f\wedge\omega_{n-1}, \hspace{2ex}  \mbox{or equivalently}\hspace{2ex} \star(\bar\partial f\wedge\omega_{n-1}) = i\bar\partial f,$$ since $\star\star = -\mbox{Id}$ on forms of odd degree.

 The contention follows.  \hfill $\Box$

 \vspace{3ex}

 \noindent {\it End of proof of Proposition \ref{Prop:commutation-defect_2-forms_integrated}.} Integrating (\ref{eqn:commutation-defect_2-forms_integrated_1}) and applying Lemma \ref{Lem:Laplacian-balanced_inner-rpod} with $f=\Lambda_\omega\alpha$, we get (\ref{eqn:commutation-defect_2-forms_integrated}). \hfill $\Box$

\vspace{3ex}

%Now, recall the following facts of [Dem84] (see also [Dem97, VII, $\S.1$]). For any Hermitian metric $\omega$ on a complex manifold $X$ with $\mbox{dim}_\C X=n$, one defines the operator $\tau=\tau_\omega:=[\Lambda_\omega,\,\partial\omega\wedge\cdot]$ of order $0$ and of type $(1,\,0)$ acting on the differential forms of $X$. The K\"ahler commutation relations generalise to the arbitrary Hermitian context as \begin{eqnarray}\label{eqn:Hermitian-commutation}i[\Lambda_\omega,\,\bar\partial] = \partial^\star + \tau^\star\end{eqnarray} and the three other relations obtained from this one by conjugation and/or adjunction. (See [Dem97, VII, $\S.1$, Theorem 1.1.].) Moreover, considering the Laplacians $$\Delta_\tau:=[d+\tau,\,d^\star + \tau^\star]  \hspace{3ex} \mbox{and} \hspace{3ex} \Delta'_\tau:=[\partial+\tau,\,\partial^\star + \tau^\star],$$ the following formula holds (see [Dem97, VII, $\S.1$, Proposition 1.16.]): \begin{eqnarray}\label{eqn:Delta_tau_sum-formula}\Delta_\tau = \Delta'_\tau + \Delta''.\end{eqnarray}

 The next consequence of the above discussion can be conveniently worded in terms of Demailly's torsion operator $\tau=\tau_\omega:=[\Lambda_\omega,\,\partial\omega\wedge\cdot]$ and the induced Laplacian $\Delta_\tau:=[d+\tau,\,d^\star + \tau^\star]$ mentioned in the introduction.

\begin{Cor}\label{Cor:Delta_tau_harmonic_implication} Let $(X,\,\omega)$ be a connected {\bf complete balanced} manifold, $\mbox{dim}_{\C}X=n\geq 2$. For any $(1,\,1)$-form $\alpha^{1,\,1}\in\mbox{Dom}\,(\Delta_\tau)$, the following implication holds: \begin{eqnarray*}\label{eqn:Delta_tau_harmonic_implication}\Delta_\tau\alpha^{1,\,1} = 0 \implies \Lambda_\omega\alpha^{1,\,1} \hspace{1ex}\mbox{is {\bf constant}}.\end{eqnarray*}

\end{Cor}

\noindent {\it Proof.} Thanks to (\ref{eqn:Delta_tau_sum-formula}) and to $\Delta'_\tau\geq 0$ and $\Delta''\geq 0$, the hypothesis $\Delta_\tau\alpha^{1,\,1} = 0$ translates to $\Delta'_\tau\alpha^{1,\,1} = 0$ and $\Delta''\alpha^{1,\,1} = 0$. Since $\omega$ is complete, these conditions are further equivalent to 
\begin{eqnarray*}\label{eqn:Delta_tau_harmonic_conditions}\nonumber &(i)& (\partial+\tau) \alpha^{1,\,1} = 0,   \hspace{15ex}  (iii)\hspace{3ex} \bar\partial\alpha^{1,\,1} = 0  \\
  &(ii)& (\partial^\star+\tau^\star) \alpha^{1,\,1} = 0,   \hspace{13ex}  (iv)\hspace{3ex} \bar\partial^\star\alpha^{1,\,1} = 0.\end{eqnarray*}

Thus, we get: \begin{eqnarray}\bar\partial(\Lambda_\omega\alpha^{1,\,1}) = [\bar\partial,\,\Lambda_\omega]\,\alpha^{1,\,1} = i(\partial^\star + \tau^\star)\,\alpha^{1,\,1} = 0,\end{eqnarray} where the first equality follows from (iii) of (\ref{eqn:Delta_tau_harmonic_conditions}), the second equality follows from Demailly's Hermitian commutation relation (\ref{eqn:Hermitian-commutation}) and the third equality follows from (ii) of (\ref{eqn:Delta_tau_harmonic_conditions}).

We conclude that the hypothesis $\Delta_\tau\alpha^{1,\,1} = 0$ implies $\bar\partial(\Lambda_\omega\alpha^{1,\,1}) = 0$. This implies, thanks to Proposition \ref{Prop:commutation-defect_2-forms_integrated} applied with $h=1$, that $\Delta(\omega_{n-1}\wedge\alpha^{1,\,1}) = 0$, where $\Delta = \Delta_\omega = dd^\star + d^\star d$ is the $d$-Laplacian induced by $\omega$. Since $\omega_{n-1}\wedge\alpha^{1,\,1} = (\Lambda_\omega\alpha^{1,\,1})\,\omega_n = \star(\Lambda_\omega\alpha^{1,\,1})$ and since $\Delta$ commutes with $\star$, we get $\Delta(\Lambda_\omega\alpha^{1,\,1})=0$. By completeness of $\omega$, this means that $d(\Lambda_\omega\alpha^{1,\,1})=0$ on $X$, hence $\Lambda_\omega\alpha^{1,\,1}$ must be constant since $X$ is connected.  \hfill $\Box$

\vspace{3ex}

An immediate consequence of Corollary \ref{Cor:Delta_tau_harmonic_implication} is that the following linear map is {\bf well defined}: \begin{eqnarray*}\label{eqn:Delta_tau_harmonic_map}T_{\omega_n}:{\cal H}_{\Delta_\tau}^{1,\,1}(X,\,\C)\longrightarrow\C, \hspace{3ex} \alpha^{1,\,1}\longmapsto\Lambda_\omega\alpha^{1,\,1}=\frac{\alpha^{1,\,1}\wedge\omega_{n-1}}{\omega_n},\end{eqnarray*} under those assumptions, where ${\cal H}_{\Delta_\tau}^{1,\,1}(X,\,\C)$ is the space of $\Delta_\tau$-harmonic $L^2_\omega$-forms of type $(1,\,1)$.

\vspace{3ex}

%\begin{Cor}\label{Cor:Delta_tau_harmonic_bal-complete-exact} Let $X$ be a compact complex {\bf balanced hyperbolic} manifold with $\mbox{dim}_\C X=n$. Let $\pi:\widetilde{X}\longrightarrow X$ be the universal cover of $X$ and $\widetilde\omega:=\pi^\star\omega$ the lift to $\widetilde{X}$ of a balanced hyperbolic metric $\omega$ on $X$.

%   There are no non-zero semi-positive $\Delta_{\widetilde\tau}$-harmonic $L^2_{\widetilde\omega}$-forms of pure type $(1,\,1)$ on $\widetilde{X}$: $$\bigg\{\alpha^{1,\,1}\in{\cal H}^{1,\,1}_{\Delta_{\widetilde\tau}}(\widetilde{X},\,\C)\,\mid\,\alpha^{1,\,1}\geq 0\bigg\} = \{0\},$$ where $\widetilde\tau=\widetilde\tau_{\widetilde\omega}:=[\Lambda_{\widetilde\omega},\,\partial\widetilde\omega\wedge\cdot]$

%\end{Cor} 

\noindent {\it Proof of Theorem \ref{The:Delta_tau_harmonic_bal-complete-exact}.} The pair $(\widetilde{X},\, \widetilde\omega)$ satisfies the hypotheses of Corollary \ref{Cor:Delta_tau_harmonic_implication} (playing the role of the pair $(X,\,\omega)$ therein). By the balanced hyperbolic hypothesis on $(X,\,\omega)$, there exists an $\widetilde\omega$-bounded smooth $(2n-3)$-form $\widetilde\Gamma$ on $\widetilde{X}$ such that $\widetilde\omega_{n-1} = d \widetilde\Gamma$.

Let $\alpha^{1,\,1}\in{\cal H}^{1,\,1}_{\Delta_{\widetilde\tau}}(\widetilde{X},\,\C)$ such that $\alpha^{1,\,1}\geq 0$. Then, $\bar\partial\alpha^{1,\,1}=0$ (by (iii) of (\ref{eqn:Delta_tau_harmonic_conditions})) and real, hence we also have $\partial\alpha^{1,\,1}=0$. Thus, $\alpha^{1,\,1}$ is $d$-closed, so \begin{eqnarray*}\label{eqn:Delta_tau_harmonic_bal-complete-exact_1}\widetilde\omega_{n-1}\wedge\alpha^{1,\,1} = d(\widetilde\Gamma\wedge\alpha^{1,\,1})\in\mbox{Im}\,d\end{eqnarray*} because $\widetilde\Gamma\wedge\alpha^{1,\,1}$ is $L^2_{\widetilde\omega}$ and $d(\widetilde\Gamma\wedge\alpha^{1,\,1})$ is again $L^2_{\widetilde\omega}$.

  On the other hand, \begin{eqnarray*}\label{eqn:Delta_tau_harmonic_bal-complete-exact_2}\widetilde\omega_{n-1}\wedge\alpha^{1,\,1} = (\Lambda_{\widetilde\omega}\alpha^{1,\,1})\,\widetilde\omega_n\in{\cal H}^{2n}_{\Delta_{\widetilde\omega}}(\widetilde{X},\,\C)\end{eqnarray*} because $\Lambda_{\widetilde\omega}\alpha^{1,\,1}$ is constant by Corollary \ref{Cor:Delta_tau_harmonic_implication}.

  Since the subspaces ${\cal H}^{2n}_{\Delta_{\widetilde\omega}}(\widetilde{X},\,\C)$ and $\mbox{Im}\,d$ of the space of $L^2_{\widetilde\omega}$-forms of degree $2n$ on $\widetilde{X}$ are orthogonal (see (d) of Theorem \ref{The:standard_complete_operators}), we deduce that $\widetilde\omega_{n-1}\wedge\alpha^{1,\,1} = 0$. Equivalently, $\Lambda_{\widetilde\omega}\alpha^{1,\,1} = 0$. This implies that $\alpha^{1,\,1}= 0$ since $\alpha^{1,\,1}\geq 0$ by hypothesis.  \hfill $\Box$

\section{Appendix}\label{section:appendix}

A key classical fact used by Gromov in [Gro91] is that some of the Lefschetz maps at the level of differential forms are {\it quasi-isometries} w.r.t. the $L^2$-inner product. We spell out the equalities involving {\it pointwise} inner products that lead to more precise statements that were used in earlier parts of our text.

Let $\omega$ be an arbitrary Hermitian metric on an arbitrary complex manifold $X$ with $\mbox{dim}_\C X=n$. As usual, for any $r=1,\dots , n$, we put $\omega_r:=\omega^r/r!$. Recall the following standard fact.

\vspace{1ex}

{\it For every $k\leq n$ and every $r\leq n-k$, the pointwise Lefschetz operator: \begin{equation}\label{eqn:pointwise-Lefschetz-map}L_\omega^r:\Lambda^kT^\star X\longrightarrow\Lambda^{k+2r}T^\star X, \hspace{3ex} L_\omega^r(\varphi)=\omega^r\wedge\varphi,\end{equation} is {\bf injective}. When $r=n-k$, $L_\omega^{n-k}$ is even {\bf bijective}.}

\vspace{1ex}

We will compare the pointwise inner products $\langle\omega_r\wedge\varphi_1,\,\omega_r\wedge\varphi_2\rangle_\omega$ and $\langle\varphi_1,\,\varphi_2\rangle_\omega$ for arbitrary $k$-forms $\varphi_1,\varphi_2\in\Lambda^kT^\star X$. We will use the following standard formula (cf. e.g. [Voi02]): \begin{equation}\label{eqn:Lr_Lambda_commutation} [L_\omega^r,\,\Lambda_\omega] = r(k-n+r-1)\,L_\omega^{r-1} \hspace{5ex} \mbox{on}\hspace{1ex}k\mbox{-forms},\end{equation} for any integer $r\geq 1$, where $\Lambda=\Lambda_\omega = (\omega\wedge\cdot)^\star$ is the adjoint of the Lefschetz operator $L_\omega$ w.r.t. the pointwise inner product $\langle\,\,,\,\,\,\rangle_\omega$ induced by $\omega$.

\vspace{2ex}

%\subsection{Case of primitive forms}\label{subsection:appendix_primitive}

\noindent {\bf (1)\, Case of primitive forms}

\vspace{2ex}

Recall that for any non-negative integer $k\leq n$, a $k$-form $\varphi$ is said to be {\it primitive w.r.t. $\omega$} (or $\omega$-{\it primitive}, or simply {\it primitive} when no confusion is likely) if it satisfies any of the following equivalent two conditions: $$\omega_{n-k+1}\wedge\varphi = 0 \iff \Lambda_\omega\varphi = 0.$$

\begin{Lem}\label{Lem:appendix_pointwise-isometry_primitive} For every $k\leq n$, every $r\leq n-k$ and any $k$-forms $\varphi_1,\varphi_2$ one of which is {\bf $\omega$-primitive}, the following identity holds: \begin{equation}\label{eqn:appendix_pointwise-isometry_primitive}\langle\omega^r\wedge\varphi_1,\,\omega^r\wedge\varphi_2\rangle_\omega = (r!)^2\,{n-k \choose r}\,\langle\varphi_1,\,\varphi_2\rangle_\omega.\end{equation}

  In particular, the analogous equality holds for the $L^2_\omega$-inner product $\langle\langle\,\,,\,\,\,\rangle\rangle_\omega$.

\end{Lem}

\noindent {\it Proof.} To make a choice, let us suppose that $\varphi_1$ is primitive. We get: \begin{eqnarray*}\langle\omega^r\wedge\varphi_1,\,\omega^r\wedge\varphi_2\rangle_\omega & = & \langle\Lambda_\omega(\omega^r\wedge\varphi_1),\,\omega^{r-1}\wedge\varphi_2\rangle_\omega \stackrel{(i)}{=} \langle[\Lambda_\omega,\,L_\omega^r]\,\varphi_1,\,\omega^{r-1}\wedge\varphi_2\rangle_\omega \\
  & \stackrel{(ii)}{=} & r(n-k-r+1)\,\langle\omega^{r-1}\wedge\varphi_1,\,\omega^{r-1}\wedge\varphi_2\rangle_\omega \\
  & \vdots & \\
  & = & r(r-1)\dots 1\, (n-k-r+1)\,(n-k-r+2)\dots (n-k)\,\langle\varphi_1,\,\varphi_2\rangle_\omega \\
  & = & r!\,\frac{(n-k)!}{(n-k-r)!}\,\langle\varphi_1,\,\varphi_2\rangle_\omega,\end{eqnarray*} where (i) follows from $\varphi_1$ being primitive, (ii) follows from the standard formula (\ref{eqn:Lr_Lambda_commutation}), the remaining equalities except for the last one follow from analogues of (i) and (ii), while the last equality proves (\ref{eqn:appendix_pointwise-isometry_primitive}). \hfill $\Box$

\vspace{2ex}

Let us also notice that, when the powers of $\omega$ are distinct, the products involved in the analogue of (\ref{eqn:appendix_pointwise-isometry_primitive}) are actually orthogonal to each other.

\begin{Lem}\label{Lem:appendix_pointwise-orthogonal_primitive} Let $r,s,k\in\N$ with $s>0$ and $k\leq n$. For any $(k-2s)$-form $u$ and any {\bf $\omega$-primitive} $k$-form $v$, the following identity holds: \begin{equation}\label{eqn:appendix_pointwise-orthogonal_primitive}\langle\omega^{r+s}\wedge u,\,\omega^r\wedge v\rangle_\omega = 0.\end{equation}

  In particular, the analogous equality holds for the $L^2_\omega$-inner product $\langle\langle\,\,,\,\,\,\rangle\rangle_\omega$.

\end{Lem}

\noindent {\it Proof.} We have: \begin{eqnarray*}\langle\omega^{r+s}\wedge u,\,\omega^r\wedge v\rangle_\omega & = & \langle\omega^{r+s-1}\wedge u,\,\Lambda_\omega(\omega^r\wedge v)\rangle_\omega \stackrel{(i)}{=} \langle\omega^{r+s-1}\wedge u,\,[\Lambda_\omega,\,L_\omega^r]\, v\rangle_\omega \\
  & \stackrel{(ii)}{=} & c_1\,\langle\omega^{r+s-1}\wedge u,\,\omega^{r-1}\wedge v\rangle_\omega = \dots = c_1\dots c_r\,\langle\omega^s\wedge u,\, v\rangle_\omega \\
%  & \vdots & \\
%  & = & c_1\dots c_r\,\langle\omega^s\wedge u,\, v\rangle_\omega \\
  & = & c_1\dots c_r\,\langle\omega^{s-1}\wedge u,\, \Lambda_\omega v\rangle_\omega = 0,\end{eqnarray*} where (i) follows from $v$ being primitive, (ii) follows from the standard formula (\ref{eqn:Lr_Lambda_commutation}) with the appropriate constant $c_1$ (whose actual value is irrelevant here), the remaining equalities except for the last one follow from analogues of (i) and (ii) with the appropriate constants $c_2, \dots , c_r$, while the last equality follows again from $v$ being primitive and proves (\ref{eqn:appendix_pointwise-orthogonal_primitive}). \hfill $\Box$

\vspace{2ex}

%\subsection{Case of primitive forms}\label{subsection:appendix_primitive}

\noindent {\bf (2)\, Case of arbitrary forms}

\vspace{2ex}

Let $\varphi_1$, $\varphi_2$ be arbitrary $k$-forms and let \begin{eqnarray}\label{eqn:appendix_Lefschetz-decomp}\varphi_1 = \varphi_{1,\,prim} + \omega\wedge\varphi_{1,\,1} + \dots + \omega^l\wedge\varphi_{1,\,l}  \hspace{2ex}\mbox{and}\hspace{2ex} \varphi_2 = \varphi_{2,\,prim} + \omega\wedge\varphi_{2,\,1} + \dots + \omega^l\wedge\varphi_{2,\,l}\end{eqnarray} be their respective {\it Lefschetz decompositions}, where $l$ is the non-negative integer defined by requiring $2l=k$ if $k$ is even and $2l=k-1$ if $k$ is odd, while the forms $\varphi_{j,\,prim}, \varphi_{j,\,1}, \dots , \varphi_{j,\,l}$ are {\it primitive} of respective degrees $k, k-2, \dots , k-2l$ for every $j\in\{1,2\}$.

\vspace{2ex}

The sense in which the Lefschetz operator (\ref{eqn:pointwise-Lefschetz-map}) is a {\it quasi-isometry} for the pointwise inner product (hence also the $L^2$-inner product) induced by $\omega$ is made explicit in the following

\begin{Lem}\label{Lem:appendix_pointwise_quasi-isometry} Fix integers $0\leq k\leq n$, $0\leq r\leq n-k$ and arbitrary $k$-forms $\varphi_1,\varphi_2$.

\vspace{1ex}

(i)\, The following identity holds: \begin{eqnarray}\label{eqn:appendix_pointwise_quasi-isometry_1}\langle\omega^r\wedge\varphi_1,\,\omega^r\wedge\varphi_2\rangle_\omega & = & (r!)^2\,{n-k \choose r}\,\langle\varphi_{1,\,prim},\,\varphi_{2,\,prim}\rangle_\omega \\
\nonumber    & + & ((r+1)!)^2\,{n-k+2 \choose r+1}\,\langle\varphi_{1,\,1},\,\varphi_{2,\,1}\rangle_\omega + \dots + ((r+l)!)^2\,{n-k+2l \choose r+l}\,\langle\varphi_{1,\,l},\,\varphi_{2,\,l}\rangle_\omega.\end{eqnarray}
%\nonumber    & \vdots & \\
%    & + & ((r+l)!)^2\,{n-k+2l \choose r+l}\,\langle\varphi_{1,\,l},\,\varphi_{2,\,l}\rangle_\omega\end{eqnarray}

%\noindent and \begin{eqnarray}\label{eqn:appendix_pointwise_quasi-isometry_2}\nonumber\langle\varphi_1,\,\varphi_2\rangle_\omega & = & \langle\varphi_{1,\,prim},\,\varphi_{2,\,prim}\rangle_\omega \\
%  & + & (1!)^2\,{n-k+2 \choose 1}\,\langle\varphi_{1,\,1},\,\varphi_{2,\,1}\rangle_\omega + \dots + (l!)^2\,{n-k+2l \choose l}\,\langle\varphi_{1,\,l},\,\varphi_{2,\,l}\rangle_\omega\end{eqnarray}
%\nonumber    & + & (1!)^2\,{n-k+2 \choose 1}\,\langle\varphi_{1,\,1},\,\varphi_{2,\,1}\rangle_\omega \\
%\nonumber    & \vdots & \\
%    & + & (l!)^2\,{n-k+2l \choose l}\,\langle\varphi_{1,\,l},\,\varphi_{2,\,l}\rangle_\omega\end{eqnarray}

\vspace{1ex}

(ii)\, Putting $C_{n,\,k,\,r,\,s}:=((r+s)!(n-k+s)!)/(s!(n-k-r+s)!)$ and $$A_{n,\,k,\,r}:=\min\limits_{s=0,\dots , l}C_{n,\,k,\,r,\,s}, \hspace{3ex} B_{n,\,k,\,r}:=\max\limits_{s=0,\dots , l}C_{n,\,k,\,r,\,s},$$ the following inequalities hold: \begin{eqnarray*}\label{eqn:appendix_pointwise_quasi-isometry_sandwich_norms}A_{n,\,k,\,r}\,|\varphi|^2_\omega\leq|\omega^r\wedge\varphi|^2_\omega\leq B_{n,\,k,\,r}\,|\varphi|^2_\omega.\end{eqnarray*}

\vspace{1ex}

(iii)\, With the notation of (ii), if $\langle\varphi_{1,\,s},\,\varphi_{2,\,s}\rangle_\omega\geq 0$ for every $s\in\{0,1,\dots , l\}$, the following inequalities hold: \begin{eqnarray*}\label{eqn:appendix_pointwise_quasi-isometry_sandwich_products}A_{n,\,k,\,r}\,\langle\varphi_1,\,\varphi_2\rangle_\omega\leq\langle\omega^r\wedge\varphi_1,\,\omega^r\wedge\varphi_2\rangle_\omega\leq B_{n,\,k,\,r}\,\langle\varphi_1,\,\varphi_2\rangle_\omega.\end{eqnarray*}

\end{Lem}

\noindent {\it Proof.} (i) Using the Lefschetz decompositions (\ref{eqn:appendix_Lefschetz-decomp}) and Lemma \ref{Lem:appendix_pointwise-orthogonal_primitive}, we get: \begin{eqnarray*}\langle\omega^r\wedge\varphi_1,\,\omega^r\wedge\varphi_2\rangle_\omega & = & \langle\omega^r\wedge\varphi_{1,\,prim},\,\omega^r\wedge\varphi_{2,\,prim}\rangle_\omega + \sum\limits_{s=1}^l\langle\omega^{r+s}\wedge\varphi_{1,\,s},\,\omega^{r+s}\wedge\varphi_{2,\,s}\rangle_\omega.\end{eqnarray*} Identity (\ref{eqn:appendix_pointwise_quasi-isometry_1}) follows from this and from Lemma \ref{Lem:appendix_pointwise-isometry_primitive}.

(ii) and (iii) follow at once from (i) applied twice, with a given $1\leq r\leq n-k$ and with $r=0$.  \hfill $\Box$

\vspace{3ex}

\noindent {\bf References.} \\

%\noindent [Bro78]\, R. Brody --- {\it Compact Manifolds and Hyperbolicity} --- Trans. Amer. Math. Soc. {\bf 235} (1978), 213-219.  

%\vspace{1ex}

\noindent [Dem84]\, J.-P. Demailly --- {\it Sur l'identit\'e de Bochner-Kodaira-Nakano en g\'eom\'etrie hermitienne} --- S\'eminaire d'analyse P. Lelong, P. Dolbeault, H. Skoda (editors) 1983/1984, Lecture Notes in Math., no. {\bf 1198}, Springer Verlag (1986), 88-97.

\vspace{1ex}

\noindent [Dem92]\, J.-P. Demailly --- {\it Regularization of Closed Positive Currents and Intersection Theory} --- J. Alg. Geom., {\bf 1} (1992), 361-409.

\vspace{1ex}

\noindent [Dem97]\, J.-P. Demailly --- {\it Complex Analytic and Algebraic Geometry} --- http://www-fourier.ujf-grenoble.fr/~demailly/books.html

\vspace{1ex}

\noindent [DGMS75]\, P. Deligne, Ph. Griffiths, J. Morgan, D. Sullivan --- {\it Real Homotopy Theory of K\"ahler Manifolds} --- Invent. Math. {\bf 29} (1975), 245-274.

\vspace{1ex}

\noindent [FLY12]\, J. Fu, J. Li, S.-T. Yau -- {\it Balanced Metrics on Non-K\"ahler Calabi-Yau Threefolds} --- J. Differential Geom. 90 (2012), p. 81–129.

\vspace{1ex}

\noindent [Fri89]\, R. Friedman --{\it On Threefolds with Trivial Canonical Bundle} --- in Complex Geometry and Lie Theory (Sundance, UT, 1989), Proc. Sympos. Pure Math., vol. 53, Amer. Math. Soc., Providence, RI, 1991, p. 103–134.

\vspace{1ex}

\noindent [Fri19]\, R. Friedman --- {\it The $\partial\bar\partial$-Lemma for General Clemens Manifolds} --- Pure and Applied Mathematics Quarterly, vol. 15, no. {\bf 4} (2019), 1001–1028.

\vspace{1ex}

\noindent [Gaf54]\, M. P. Gaffney --- {\it A Special Stokes's Theorem for Complete Riemannian Manifolds} --- Ann. of Math. {\bf 60}, No. 1, (1954), 140-145.

\vspace{1ex}

\noindent [Gau77]\, P. Gauduchon --- {\it Le th\'eor\`eme de l'excentricit\'e nulle} --- C. R. Acad. Sci. Paris, S\'er. A, {\bf 285} (1977), 387-390.

%\vspace{1ex}

%\noindent [Gau77b]\, P. Gauduchon --- {\it Fibr\'es hermitiens \`a endomorphisme de Ricci non n\'egatif} --- Bull. Soc. Math. France {\bf 105} (1977) 113-140.

\vspace{1ex}

\noindent [Gro91]\, M. Gromov --- {\it K\"ahler Hyperbolicity and $L^2$ Hodge Theory} --- J. Diff. Geom. {\bf 33} (1991), 263-292.

%\vspace{1ex}

%\noindent [Kob70]\, S. Kobayashi --- {\it Hyperbolic Manifolds and Holomorphic Mappings} --- Marcel Dekker, New York (1970). 

\noindent [Lam99]\, A. Lamari --- {\it Courants k\"ahl\'eriens et surfaces compactes} --- Ann. Inst. Fourier {\bf 49}, no. 1 (1999), 263-285.

%\vspace{1ex}

%\noindent [Lan87]\, S. Lang --- {\it Introduction to Complex Hyperbolic Spaces} --- Springer-Verlag (1987).

\vspace{1ex}

\noindent [LT93]\, P. Lu, G. Tian -- {\it The Complex Structures on Connected Sums of $S^3\times S^3$} --- in Manifolds and Geometry (Pisa, 1993), Sympos. Math., XXXVI, Cambridge Univ. Press, Cambridge, 1996, p. 284–293.

\vspace{1ex}

\noindent [MP21]\, S. Marouani, D. Popovici --- {\it Balanced Hyperbolic and Divisorially Hyperbolic Compact Complex Manifolds} --- arXiv e-print math.CV/2107.08972v1.

\vspace{1ex}

\noindent [Mic83]\, M. L. Michelsohn --- {\it On the Existence of Special Metrics in Complex Geometry} --- Acta Math. {\bf 143} (1983) 261-295.

\vspace{1ex}

\noindent [Pop14]\, D. Popovici --- {\it Deformation Openness and Closedness of Various Classes of Compact Complex Manifolds; Examples} --- Ann. Sc. Norm. Super. Pisa Cl. Sci. (5), Vol. XIII (2014), 255-305.

\vspace{1ex}

\noindent [Pop15]\, D. Popovici --- {\it Aeppli Cohomology Classes Associated with Gauduchon Metrics on Compact Complex Manifolds} --- Bull. Soc. Math. France {\bf 143}, no. 4 (2015), p. 763-800.

\vspace{1ex}

\noindent [Pop16] \, D. Popovici --- {\it Sufficient Bigness Criterion for Differences of Two Nef Classes} --- Math. Ann. {\bf 364} (2016), 649-655.

%\vspace{1ex}

%\noindent [PU18]\, D. Popovici, L. Ugarte --- {\it Compact Complex Manifolds with Small Gauduchon Cone} --- Proc. LMS {\bf 116}, no. 5 (2018) doi:10.1112/plms.12110.

\vspace{1ex}

\noindent [Voi02]\, C. Voisin --- {\it Hodge Theory and Complex Algebraic Geometry. I.} --- Cambridge Studies in Advanced Mathematics, 76, Cambridge University Press, Cambridge, 2002.

\vspace{1ex}

\noindent [Yac98]\, A. Yachou --- {\it Sur les vari\'et\'es semi-k\"ahl\'eriennes} --- PhD Thesis (1998), University of Lille.

\vspace{3ex}

\noindent Institut de Math\'ematiques de Toulouse, Universit\'e Paul Sabatier,

\noindent 118 route de Narbonne, 31062 Toulouse, France

\noindent Email: almarouanisamir@gmail.com AND popovici@math.univ-toulouse.fr

\vspace{2ex}

\noindent {\it For the first-named author only, also}:

\vspace{1ex}

\noindent Universit\'e de Monastir, Facult\'e des Sciences de Monastir

\noindent Laboratoire de recherche Analyse, G\'eom\'etrie et Applications LR/18/ES/16

\noindent Avenue de l'environnement 5019, Monastir, Tunisie

\end{document}